\documentclass[11pt]{amsart}

\usepackage[utf8]{inputenc}
\usepackage[english]{babel}
\usepackage[a4paper]{geometry}
\usepackage{pdflscape}
\usepackage{rotating}

\usepackage{amsfonts}
\usepackage{amssymb}
\usepackage{amsmath}

\usepackage[arrow,curve,matrix,tips,frame]{xy}

\usepackage{adjustbox}

\usepackage[dvipsnames,table]{xcolor}
\usepackage[most]{tcolorbox}
\usepackage{textcomp}
\usepackage{multirow}

\usepackage{fancyvrb}

\usepackage{newunicodechar} \usepackage{pmboxdraw} \usepackage{pifont} \newunicodechar{∩}{\ensuremath{\cap}} \newunicodechar{∪}{\ensuremath{\cup}} \newunicodechar{⊂}{\ensuremath{\subset}} \newunicodechar{⊆}{\ensuremath{\subseteq}} \newunicodechar{≤}{\ensuremath{\leq}} \newunicodechar{≥}{\ensuremath{\geq}} \newunicodechar{ℙ}{\ensuremath{\mathbb P}} \newunicodechar{μ}{\ensuremath{\mu}} \newunicodechar{χ}{\ensuremath{\chi}} \newunicodechar{Ũ}{\ensuremath{\widetilde U}} \newunicodechar{→}{\ensuremath{\to}} \newunicodechar{⟶}{\ensuremath{\longrightarrow}} \newunicodechar{⇢}{\ensuremath{\dashrightarrow}} \newunicodechar{⁻}{\ensuremath{{}^{-}}} \newunicodechar{¹}{\ensuremath{{}^{1}}} \newunicodechar{²}{\ensuremath{{}^{2}}} \newunicodechar{³}{\ensuremath{{}^{3}}} \newunicodechar{⁴}{\ensuremath{{}^{4}}} \newunicodechar{⁵}{\ensuremath{{}^{5}}} \newunicodechar{⁶}{\ensuremath{{}^{6}}} \newunicodechar{⁷}{\ensuremath{{}^{7}}} \newunicodechar{⁸}{\ensuremath{{}^{8}}} \newunicodechar{⁹}{\ensuremath{{}^{9}}} \newunicodechar{✓}{\ding{51}} \newunicodechar{✅}{\ding{51}} \newunicodechar{★}{\ensuremath{\star}} \newunicodechar{✦}{\ensuremath{\diamond}}
\newunicodechar{≃}{\ensuremath{\simeq}} \newunicodechar{≅}{\ensuremath{\cong}}

\definecolor{airforceblue}{rgb}{0.2,0.2,0.6}
\definecolor{bleudefrance}{rgb}{0.0,0.0,1.0}
\definecolor{darkorange}{rgb}{1.0,0.55,0.0}
\definecolor{commentoutput}{rgb}{0.33,0.38,0.44}
\definecolor{output}{rgb}{0.24,0.08,0.08}
\definecolor{circOut}{rgb}{0.0,0.81,0.82}

\def\gt;{>}
\def\lt;{<}
\def\amp;{\char38}

\DefineVerbatimEnvironment{M2Verbatim}{Verbatim}{ frame=single, framesep=2mm, rulecolor=\color{gray!55}, fillcolor=\color{gray!4} }

\usepackage[backref=page]{hyperref}

\hypersetup{
  pdfborder={0 0 0},
  colorlinks=true,
  linktoc=page,
  pdfauthor={F. Russo and G. Staglian{\`o}},
  pdftitle={Rational cubic fourfolds containing a plane}
}

\newtheorem{thm}{Theorem}[section]
\newtheorem{cor}[thm]{Corollary}
\newtheorem{Proposition}[thm]{Proposition}
\newtheorem{defin}[thm]{Definition}
\newtheorem{Lemma}[thm]{Lemma}
\newtheorem{rmk}[thm]{Remark}
\newtheorem*{mainthm}{Main Theorem}

\theoremstyle{definition}
\newtheorem{ex}[thm]{Example}

\def\p{\mathbb P}
\def\P{\mathbb P}
\def\H{\operatorname{H}}
\def\rk{\operatorname{rk}}
\def\h{\operatorname{h}}
\def\Sing{\operatorname{Sing}}
\def\mult{\operatorname{mult}}
\def\Rat{\operatorname{Rat}}
\def\Pic{\operatorname{Pic}}
\def\codim{\operatorname{codim}}
\def\Bl{\operatorname{Bl}}
\def\disc{\operatorname{disc}}
\def\Bs{\operatorname{Bs}}

\newcommand{\PGL}[1]{\operatorname{PGL}_{#1}}

\numberwithin{equation}{section}

\begin{document}

\title[Rational cubic fourfolds containing a plane]{Some rational cubic fourfolds containing a plane and their associated K3 surfaces}

\subjclass[2020]{Primary 14E08; Secondary 14M20, 14M07, 14N05, 14J28, 14J70}
\keywords{Rationality of cubic fourfolds, special cubic fourfolds, K3 surfaces, moduli spaces}

\author[F. Russo]{Francesco Russo}
\address{Dipartimento di Matematica e Informatica, Universit\`a degli Studi di Catania,
Viale A. Doria 5, 95125 Catania, Italy}
\email{francesco.russo@unict.it}
\thanks{The first author was partially supported by the NextGenEU under the MUR PRIN 2022 program (Grant No. 20223B5S8L, project title: ``Birational geometry of moduli spaces and special varieties'', CUP E53D23005410006).}

\author[G. Staglian\`o]{Giovanni Staglian\`o}
\address{Dipartimento di Matematica e Informatica, Universit\`a degli Studi di Catania,
Viale A. Doria 5, 95125 Catania, Italy}
\email{giovanni.stagliano@unict.it}
\thanks{The second author is partially supported by the PRIN 2022 project ``Unirationality, Hilbert schemes, and singularities'' (MUR, DD n.~104/2022 and DD n.~1401/2024; CUP E53C24002310006).}
\thanks{Both authors are members of the G.N.S.A.G.A. of INdAM\@.}

\begin{abstract}
We give geometric descriptions of general members of the irreducible divisors
$\mathcal S_\delta\subset\mathcal C_8$ for
$\delta=21,29,37,45,53,61,69$,
where $\mathcal C_8$ is the divisor in the moduli space $\mathcal C$ of cubic fourfolds parametrizing those containing a plane, and we prove that these divisors are unirational. Sections of the quadric fibration induced by projection from the plane give rise to birational maps from these cubic fourfolds to projective models of $\p^2\times\p^2$. The same maps yield explicit projective models of the K3 surfaces in the corresponding divisors $\mathcal K_\delta$ of the moduli space $\mathcal K$ of polarized K3 surfaces of degree $2$; a very general member of $\mathcal K_\delta$ has Picard rank $2$ and Picard lattice of discriminant $-\delta$.
\end{abstract}

\maketitle

\section*{Introduction}

The longstanding problem of the irrationality of a very general cubic fourfold has recently been settled in \cite{KKPY} using the theory of Hodge atoms. The locus $\Rat(\mathcal C)$ of rational cubic fourfolds in the quasi-projective moduli space $\mathcal C$ is a countable union of irreducible closed subsets; see \cite{dFF,KontsevichTschinkelInventiones}. The irreducible divisor $\mathcal C_8\subset\mathcal C$ parametrizes cubic fourfolds $X$ containing a plane $P$, which is unique when $X$ is general in $\mathcal C_8$. Moreover, a very general cubic fourfold in $\mathcal C_8$ is irrational; see \cite[Example 6.17]{KKPY}. The methods developed in \cite{KKPY}, however, do not describe $\Rat(\mathcal C)$ or $\Rat(\mathcal C)\cap\mathcal C_8$ explicitly, beyond locating them in the corresponding Noether--Lefschetz loci.

Let $Q$ be a general fibre of the quadric fibration obtained by restricting to $X$ the projection from $P$. Hassett proved in \cite{HassettPlane} that, for infinitely many integers $\delta>13$ satisfying $\delta\equiv 5\pmod 8$, there exists an irreducible divisor $\mathcal S_\delta\subset\mathcal C_8$ whose very general member contains a $2$-cycle $T$ such that $T\cdot Q$ is odd. Here, $\delta$ is the discriminant of the intersection matrix of the rank-three lattice $\langle h^2,Q,T\rangle$, where $h$ denotes the class of a hyperplane section. It is expected that every $\delta\geq 21$ satisfying this congruence condition occurs; see \cite[Remark 4.3, 2)]{HassettPlane}. These results are obtained through the analysis of the period map and infinitesimal Hodge theory, methods also used in \cite{CHM,Voisin2,AHTVA,HPT2}. The cycle $T$ induces a $0$-cycle of odd degree on the generic fibre of the quadric fibration. This implies the existence of a rational section and hence the rationality of a general cubic fourfold in $\mathcal S_\delta$. Therefore,
\[
\mathcal S_\delta\subset\Rat(\mathcal C)\cap\mathcal C_8
\]
by \cite{KontsevichTschinkelInventiones}. It is also known that $\mathcal C_8\cap\Rat(\mathcal C)$ contains irreducible closed subvarieties of dimension $18$ that are not among the divisors $\mathcal S_\delta$; see, for example, \cite[Remark 4.3, 1)]{HassettPlane} and \cite{BRS}.

Our approach to the rationality of these cubic fourfolds is geometric. If $S\subset X$ is a rational section of the quadric fibration, then $P\cup S$ is a reducible surface with one apparent double point, or reducible OADP surface: through a general point of $\p^5$ there passes a unique secant line to $P\cup S$. A general secant line meets $X$ in exactly one further point, so that $X$ is birational to the variety of secant lines to $P\cup S$, which is in turn birational to $P\times S$, and hence to $\p^2\times\p^2$. This is analogous to the classical treatment of the Fano--Morin divisor, also known as $\mathcal C_{14}$, where the rationality of every cubic fourfold is proved using reducible or irreducible OADP surfaces; see \cite{Fano,BD,BRS}. Congruences of secant curves to suitable surfaces have also been used to prove the rationality of the cubic fourfolds in other Hassett divisors and to construct the associated K3 surfaces explicitly; see \cite{RS1,Explicit,JEMS,S-expMath}. Reducible and irreducible OADP varieties are also useful over non-algebraically closed fields: for example, a smooth cubic surface is rational over its base field if it contains either a pair of skew lines or a twisted cubic defined over that field.

It is not difficult to construct families of reducible OADP surfaces in $\p^5$ of the form $P\cup S$, or equivalently sections of the projection of $\p^5$ from $P$; see Proposition \ref{degpiP}. However, only a few members of these families are contained in cubic fourfolds, as in the case of reducible OADP curves in $\p^3$ discussed in Example \ref{congP3}. The results of \cite{HassettPlane}, together with Proposition \ref{degpiP}, ensure the existence of countably many distinct families of reducible OADP surfaces $P\cup S_\delta$ contained in general cubic fourfolds belonging to $\mathcal S_\delta$, for infinitely many integers $\delta>13$ satisfying $\delta\equiv 5\pmod 8$. An explicit description of the entire series of surfaces $S_\delta$ is still unknown.

In this paper we construct explicit sections \( S\subset X\subset\p^5 \) of the quadric fibrations associated with general points $[X]\in\mathcal S_\delta$, for \( \delta=21,29,37,45,53,61,69, \) and we describe the corresponding divisors $\mathcal K_\delta$ of associated K3 surfaces in the moduli space of polarized K3 surfaces of degree $2$. Our discussion begins with the two classical cases described in \cite{HassettPlane}: the divisor $\mathcal S_{21}$, whose general member contains two disjoint planes, and the divisor $\mathcal S_{45}$, studied in \cite{Tregub}, whose general member contains a Veronese surface together with one of its trisecant planes. A selection of examples is collected in Table \ref{TABex}, and our main results can be summarized as follows.

\begin{mainthm}
For each $\delta\in\{21,29,37,45,53,61,69\}$, the following hold.
\begin{enumerate}
\item[{\rm (1)}] A general cubic fourfold $[X]\in\mathcal S_\delta$ contains a reducible OADP surface $P\cup S$, where $P$ is a plane and $S$ is a rational surface with at most nodes as singularities, such that the pair $(S,P)$ is of one of the types listed in Table \ref{TABex}. In particular, $S$ is a section of the quadric fibration of $X$ induced by projection from $P$.

\item[{\rm (2)}] The divisor $\mathcal S_\delta$ is unirational, and so is the corresponding divisor $\mathcal K_\delta$ of associated K3 surfaces in the moduli space $\mathcal K$ of polarized K3 surfaces of degree $2$.

\item[{\rm (3)}] For a general $[X]\in\mathcal S_\delta$, the secant lines to $P\cup S$ give an explicit birational map from $X$ to a projective model of $\p^2\times\p^2$ and an explicit smooth projective model $\widetilde U$ of the associated K3 surface $K_X$. The hyperplane class $H_X$ of this model and the degree-$2$ polarization $D_X$ of $K_X$ span a lattice of discriminant $-\delta$.
\end{enumerate}
\end{mainthm}
The proof combines a general geometric construction with a parameter count and it may be extended to cover also other few cases, e.g. $\delta=77$ (see Remark \ref{refinedCount}). By Proposition \ref{degpiP}, sections of the quadric fibration correspond to reducible OADP surfaces $P\cup S$, and the construction of Example \ref{Yer} produces such surfaces from plane Cremona transformations; their discriminant $\delta$ is given by formula \eqref{deltaex}. The parameter count of Section \ref{ParCount}, and in particular Proposition \ref{suffSd}, shows that the cubic fourfolds containing a general surface $P\cup S$ of a given type fill the divisor $\mathcal S_\delta$. This yields (1) and the unirationality of $\mathcal S_\delta$, and hence that of its image $\mathcal K_\delta$ under the rational map $[X]\mapsto[(K_X,D_X)]$. In Section \ref{SecantLines}, the K3 surface $K_X$ is recovered from the surface $U_S$ parametrizing the secant lines to $P\cup S$ contained in $X$: Proposition \ref{doublecoverT} describes $U_S$ in terms of a blow-up of $K_X$, and the Fano maps of $P\cup S$, which contract the secant lines, produce the models $\widetilde U$ in (3) through Proposition \ref{US=U}. The facts that depend on the specific type of $P\cup S$, such as the existence and smoothness of cubic fourfolds containing $P\cup S$ and the values of $\h^0(\mathcal I_{S\cup P}(3))$ and $\h^0(N_{S/X})$, are verified with \emph{Macaulay2} \cite{macaulay2}, using the package \emph{SpecialFanoFourfolds} \cite{SpecialFanoFourfoldsSource,2024packageSFF}, which also computes the invariants of $W$, $U$ and $\widetilde U$; see Table \ref{TABex} and Section \ref{Computations}.

The explicit description of the divisors $\mathcal S_\delta$ is part of the broader problem of describing $\Rat(\mathcal C)$. Each $\mathcal S_\delta$ is a codimension-two locus of rational cubic fourfolds in $\mathcal C$, and it is contained in infinitely many Hassett divisors $\mathcal C_d$; see Section \ref{SdeltaSec}. In some of our examples, the class of $S$ has discriminant $\disc\langle h^2,S\rangle=62$, $78$ or $86$; for instance, the unirational divisor $\mathcal S_{53}$ is an irreducible component of $\mathcal C_8\cap\mathcal C_{86}$. See Section \ref{S53}, and \cite{FL} for other examples of this kind. A major obstacle to a geometric description of the sections $S_\delta$ for large $\delta\equiv 5\pmod 8$, as well as to an explicit construction of the corresponding divisors in moduli spaces of polarized K3 surfaces, may be the failure of $\mathcal S_\delta$ to be uniruled. Recall that $\mathcal C_d$ is unirational for several small values of $d$; see \cite{Nuer-Unir} and the references therein. On the other hand, $\mathcal C_d$ has nonnegative Kodaira dimension, and hence is not uniruled, for $d=6n+2$ with $n>13$ and $n\neq 15$ \cite{TaVaA}; this includes every admissible $d\geq 86$ with $d\equiv 2\pmod 6$.

The paper is organized as follows. In Section \ref{Preliminaries}, we recall some facts on Hassett divisors, on sections of projections from linear spaces, and on the divisors $\mathcal S_\delta$. Section \ref{Examples} contains some general facts on reducible OADP varieties, and in Section \ref{Exseries} we construct the families of reducible OADP surfaces $P\cup S$ used throughout the paper. In Section \ref{ParCount}, we compute the dimension of the corresponding loci of cubic fourfolds, which proves the first two parts of the Main Theorem. Section \ref{SecantLines} is devoted to the secant lines to $P\cup S$ contained in $X$ and to the construction of the associated K3 surfaces. Finally, in Section \ref{Computations} we illustrate some of the functionalities of the \emph{Macaulay2} package \emph{SpecialFanoFourfolds} \cite{SpecialFanoFourfoldsSource,2024packageSFF}, developed by the second-named author, which were used for the computations in this paper; Section \ref{SummaryTables} contains the tables of examples.

\section{Preliminaries}\label{Preliminaries}

We recall some definitions and known facts.

\subsection{Hassett divisors and the self-intersection formula}

Let $\mathcal C$ denote the $20$-dimensional quasi-projective moduli space of smooth cubic fourfolds in $\p^5$. For a cubic fourfold $X$, let $h$ denote the class of a hyperplane section and set
\[
\H^{2,2}(X,\mathbb Z)
=
\H^{2,2}(X)\cap\H^4(X,\mathbb Z).
\]
For every positive integer $d$, the corresponding Noether--Lefschetz locus is
\[
\mathcal C_d
=
\left\{
[X]\in\mathcal C:
\begin{array}{c}
\text{there exists a primitive rank-two lattice }
K\subseteq\H^{2,2}(X,\mathbb Z)\\
\text{containing }h^2\text{ and having discriminant }d
\end{array}
\right\}.
\]
If $K=\langle h^2,S\rangle$ for an algebraic $2$-cycle $S$ on $X$, then
\[
\disc(K)
=
\det
\begin{pmatrix}
3 & h^2\cdot S\\
h^2\cdot S & S^2
\end{pmatrix}
=
3S^2-(h^2\cdot S)^2.
\]

Hassett proved that each nonempty locus $\mathcal C_d$ is an irreducible divisor, and that
\[
\mathcal C_d\neq\emptyset
\quad\Longleftrightarrow\quad
d>6
\quad\text{and}\quad
d\equiv 0,2\pmod 6.
\]
Moreover, $\H^{2,2}(X,\mathbb Z)$ has rank two for a very general $[X]\in\mathcal C_d$. Such an integer $d$ is called \emph{admissible} if it is not divisible by $4$, $9$ or any odd prime $p\equiv 2\pmod 3$; these are the values of $d$ for which the cubic fourfolds in $\mathcal C_d$ have associated K3 surfaces in the sense of Hodge theory. See \cite{Has00} and the survey \cite{Levico} for these results.

Let $S\subset\p^5$ be an irreducible surface with $\nu\geq 0$ nodes, and let $\widetilde S$ be its smooth normalization. Suppose that a general hyperplane section of $S$ is a smooth irreducible curve of genus $g(S)$, and let $X\subset\p^5$ be a cubic fourfold containing $S$. The self-intersection of $S$ in $X$ is then given by
\begin{equation}\label{selfS}
S^2
=
3\deg(S)+6g(S)-6
+2K_{\widetilde S}^2
-12\chi(\mathcal O_{\widetilde S})
+2\nu.
\end{equation}
If the class of $S$ in $\H^{2,2}(X,\mathbb Z)$ is not a multiple of $h^2$, then \( \disc\langle h^2,S\rangle = 3S^2-\deg(S)^2. \) Cubic fourfolds belonging to some $\mathcal C_d$ are called \emph{special}; equivalently, they contain an algebraic surface whose class is not a multiple of $h^2$. For small values of $d$, the divisor $\mathcal C_d$ can often be described geometrically as the locus of cubic fourfolds containing an irreducible surface $S\subset\p^5$ with prescribed invariants. Such a description need not be unique.

For example, $\mathcal C_8$ is the locus of cubic fourfolds containing a plane $P$; in this case,
\[
P^2=3
\qquad\text{and}\qquad
h^2\cdot P=\deg(P)=1,
\]
so that
\[
\disc\langle h^2,P\rangle=3P^2-1=8.
\]
Every cubic fourfold in $\mathcal C_8$ contains a plane, which is unique for a general member of $\mathcal C_8$.

\subsection{Sections of projections and reducible OADP varieties}

Let $P=\p^r\subsetneq\p^N$ be a linear subspace, and let
$M=\p^{N-r-1}\subset\p^N$ be a complementary linear subspace, so that
$\langle P,M\rangle=\p^N$. We denote by \( \pi_P:\p^N\dashrightarrow M \) the projection from $P$ onto $M$.

A rational map \( s:M\dashrightarrow\p^N \) is a rational section of $\pi_P$ if $\pi_P\circ s$ is the identity on $M$. In this case, the closure
\[
Y=\overline{s(M)}\subset\p^N
\]
is a rational variety of dimension $N-r-1$, and we shall simply say that $Y$ is a section of $\pi_P$.

Let $Z^n,Y^m\subset\p^N$ be irreducible projective varieties. Their \emph{join} is the irreducible variety
\[
J(Z,Y)
=
\overline{
\bigcup_{\substack{z\in Z,\ y\in Y\\ z\neq y}}
\langle z,y\rangle
}
\subseteq\p^N
\]
swept out by the lines joining points of $Z$ and $Y$. One has \( \dim J(Z,Y)\leq n+m+1, \) and the dimension of the join is determined by the relative position of the tangent spaces to $Z$ and $Y$ at general points, according to Terracini's lemma.

Suppose that $N=n+m+1$. We denote by $\alpha_{Z,Y}$ the number of lines joining $Z$ and $Y$ and passing through a general point of $\p^N$, setting $\alpha_{Z,Y}=0$ when $J(Z,Y)\subsetneq\p^N$.

Recall that OADP stands for \emph{one apparent double point}: the number of apparent double points of a variety $Z^n\subset\p^{2n+1}$ is the number of double points acquired by its projection from a general point to $\p^{2n}$. A reducible variety
\[
Z^n\cup Y^m\subset\p^{n+m+1}
\]
is called a \emph{reducible OADP variety} if $Z\neq Y$ and $\alpha_{Z,Y}=1$. An irreducible variety \( Z^n\subset\p^{2n+1} \) is called an \emph{OADP variety} if $J(Z,Z)=\p^{2n+1}$ and $\alpha_{Z,Z}=1$. For varieties of dimensions one and two, we use the terms OADP curve and OADP surface, respectively.

If $P=\p^r\subset\p^N$ and $Y\subset\p^N$ is an irreducible variety of dimension $N-r-1$, then \( J(P,Y)=\p^N \) if and only if the restriction \( \pi_P|_Y:Y\dashrightarrow M \) is dominant.

As usual, we set the degree of a rational map between varieties of the same dimension equal to zero when the map is not dominant.

The following elementary result is likely well known, but we include it because it will be combined with the main result of \cite{HassettPlane}.

\begin{Proposition}\label{degpiP}
With the notation above, assume that $r\geq 1$. Then \( \deg(\pi_P|_Y)=\alpha_{P,Y}. \) In particular, an irreducible closed subvariety \( Y^{N-r-1}\subset\p^N \) is a rational section of the projection from $P=\p^r$ if and only if $P\cup Y$ is a reducible OADP variety.

\end{Proposition}

\begin{proof}
The case $\alpha_{P,Y}=0$ is immediate. Assume therefore that $\pi_P|_Y$ is dominant, and let \( k=\deg(\pi_P|_Y). \) For a general point $q\in\p^N$, the fibre of $\pi_P|_Y$ over $\pi_P(q)$ consists of $k$ points
\[
\langle q,P\rangle\cap Y\setminus(P\cap Y)
=
\{y_1,\ldots,y_k\}.
\]
For each $j$, the line $\langle q,y_j\rangle$ meets $P$ at a point distinct from $y_j$. Thus, each point $y_j$ determines a line joining $P$ and $Y$ and passing through $q$, so that \( k\leq\alpha_{P,Y}. \) Conversely, every line joining $P$ and $Y$ and passing through $q$ is contained in $\langle q,P\rangle$. Its point of intersection with $Y$ therefore belongs to
\[
\langle q,P\rangle\cap Y\setminus(P\cap Y).
\]
Hence \( \alpha_{P,Y}\leq k, \) and consequently
\[
\deg(\pi_P|_Y)=k=\alpha_{P,Y}.
\]

It follows that $Y$ is a rational section of $\pi_P$ if and only if \( \deg(\pi_P|_Y)=1, \) which, by the first part, is equivalent to $\alpha_{P,Y}=1$. This is precisely the condition that $P\cup Y$ be a reducible OADP variety.

\end{proof}

\subsection{The divisors $\mathcal S_\delta\subset\mathcal C_8$}\label{SdeltaSec}

Let $X\subset\p^5$ be a cubic fourfold containing a plane $P$, and let $S\subset X$ be a section of the restriction \( \pi_P|_X:X\dashrightarrow\p^2. \) The lattice $\langle h^2,P,S\rangle$ has rank three, and we define
\begin{equation}\label{defdelta}
\delta(S)
=
\disc\langle h^2,P,S\rangle
=
\det
\begin{pmatrix}
3 & 1 & \deg(S)\\
1 & 3 & P\cdot S\\
\deg(S) & P\cdot S & S^2
\end{pmatrix}.
\end{equation}

\begin{rmk}\label{TvsY}{\rm 
Let $Q\subset X$ be a fibre of the quadric fibration induced by projection from $P$. Since \( h^2=P+Q, \) we have \( \langle h^2,P\rangle=\langle h^2,Q\rangle. \) The conics
\[
Q\cap P\subset P
\]
form a two-dimensional base-point-free family. Indeed, after blowing up $P$, the quadrics are the fibres of the morphism \( f_P:\Bl_PX\longrightarrow\P^2. \) If the family of conics had a base point $q\in P$, then all the fibres of $f_P$ would meet the line over $q$ in the exceptional divisor. Infinitely many fibres would therefore pass through a general point of this line, which is impossible.

For every section $S\subset X$ of $\pi_P|_X$, the intersection number $Q\cdot S$ is odd. More precisely, if $P\cap S$ is zero-dimensional, then \( Q\cdot S=1. \) If $P\cap S$ contains a curve of degree $e\geq 1$, then Bézout's theorem gives \( Q\cdot S=2e+1. \) Conversely, suppose that $X$ contains an algebraic $2$-cycle $T$ such that $T\cdot Q$ is odd. Then $T$ induces a $0$-cycle of odd degree on the generic fibre of the quadric fibration. By \cite[Proposition 2.1]{HassettPlane}, the fibration admits a rational section.

Thus, if the quadric fibration associated with a cubic fourfold $[X]\in\mathcal C_8$ admits a section $S$, then $\H^{2,2}(X,\mathbb Z)$ contains the rank-three lattice \( \langle h^2,Q,S\rangle \) with $Q\cdot S$ odd. Conversely, if $\H^{2,2}(X,\mathbb Z)$ contains a rank-three lattice \( \langle h^2,Q,T\rangle \) with $Q\cdot T$ odd, then the quadric fibration admits a section $S$. By Proposition \ref{degpiP}, the cubic fourfold $X$ then contains the reducible OADP surface $P\cup S$.
}
\end{rmk}

This leads to the following loci, introduced by Hassett:
\[
\mathcal{S}_\delta=\left\{[X]\in\mathcal{C}_8 :
\begin{array}{l}
\text{there exists a saturated  rank-three sublattice}\\
K=\langle h^2,Q,T\rangle\subseteq H^{2,2}(X,\mathbb{Z})
\text{ with } Q\cdot T \text{ odd and } \operatorname{disc}(K)=\delta
\end{array}
\right\}.
\]
The results of \cite{HassettPlane} show that every nonempty $\mathcal S_\delta$ is an irreducible divisor in $\mathcal C_8$, and that
\[
\mathcal S_\delta\neq\emptyset
\quad\Longrightarrow\quad
\delta\geq 21
\quad\text{and}\quad
\delta\equiv 5\pmod 8.
\]
For a very general $[X]\in\mathcal S_\delta$, one has \( \H^{2,2}(X,\mathbb Z)=\langle h^2,P,S\rangle, \) where $S\subset X$ is a section of the quadric fibration induced by $P$. Moreover, $\mathcal S_\delta$ is nonempty for infinitely many integers $\delta$ satisfying the above conditions.

It is expected that $\mathcal S_\delta$ is nonempty for every integer
\[
\delta\geq 21,
\qquad
\delta\equiv 5\pmod 8;
\]
see \cite[Remark 4.3, 2)]{HassettPlane}.

Since $\H^{2,2}(X,\mathbb Z)$ has rank three for a very general $[X]\in\mathcal S_\delta$, the divisor $\mathcal S_\delta$ is contained in infinitely many Hassett divisors. Indeed, if $x$ and $y$ are coprime integers, then the lattice $\langle h^2,xP+yS\rangle$ is primitive, and
\[
\disc\langle h^2,xP+yS\rangle
=
8x^2+2\bigl(3\,P\cdot S-\deg(S)\bigr)xy+\bigl(3S^2-\deg(S)^2\bigr)y^2.
\]
This is a positive definite binary quadratic form of discriminant $-12\delta$, so it takes infinitely many values on pairs of coprime integers. For each such value $d\neq 8$, the divisor $\mathcal S_\delta$ is an irreducible component of $\mathcal C_8\cap\mathcal C_d$, because it has dimension $18$ and $\mathcal C_8\neq\mathcal C_d$. For instance, $\mathcal S_{53}$ is contained in $\mathcal C_d$ for $d=20,26,30,48,56,84,86,\ldots$; see Section \ref{S53}.

\section{Reducible OADP varieties}\label{Examples}

We collect some general facts on reducible OADP varieties, which show some of the geometric subtleties of the problem.

\begin{rmk}\label{intOADP}{\rm
Suppose that
\[
Z^n\cup Y^m\subset\p^{n+m+1}
\]
is a reducible variety, with $n,m\geq 1$ and $N=n+m+1$. Let
\[
a=\deg(Z)
\qquad\text{and}\qquad
b=\deg(Y).
\]
For a general point $q\in\p^N$, consider the projection \( \pi_q:\p^N\dashrightarrow\p^{N-1}. \) Then
\[
\alpha_{Z,Y}
=
ab-
\sum_{C\subseteq\pi_q(Z\cap Y)}
\mult_C\bigl(\pi_q(Z)\cdot\pi_q(Y)\bigr),
\]
where the intersection product is computed in $\p^{N-1}$ and the sum runs over the irreducible components $C$ of $\pi_q(Z\cap Y)$.

Consequently, $Z\cup Y$ is a reducible OADP variety if and only if the contribution of $\pi_q(Z\cap Y)$ to the intersection product \( \pi_q(Z)\cdot\pi_q(Y) \) has degree $ab-1$, counted with multiplicity. If $Z\cap Y$ is zero-dimensional, this simply means that it is a $0$-cycle of degree $ab-1$. However, $Z\cap Y$ may also have irreducible components of positive dimension; see Remark \ref{N34} and the examples below.
}
\end{rmk}

We also use the following elementary fact.

\begin{Lemma}\label{linearX}
Let $Z\cup Y\subset\p^N$ be a reducible OADP variety whose components have positive dimension, and suppose that
\[
\dim(Z\cap Y)\leq 0.
\]
If $a=\deg(Z)$ and $b=\deg(Y)$, then one of the two varieties is a linear space, and its intersection with the other has length $ab-1$.
\end{Lemma}

\begin{proof}
If $Z\cap Y=\emptyset$, then Remark \ref{intOADP} gives $ab=1$, and hence $a=b=1$. We may thus assume that $Z\cap Y$ is a nonempty $0$-dimensional scheme of length $ab-1$.

If $a=1$, then $Z$ is a linear space and there is nothing to prove. Suppose that $a\geq 2$, and let $z\in Z$ be a general point. Projection from $z$ maps $Z$ onto a variety of degree $a-1$, while the image of $Y$ has degree $b$. These two images have complementary dimensions in $\p^{N-1}$ and therefore intersect in $(a-1)b$ points, counted with multiplicity.

On the other hand, their intersection contains the images of the $ab-1$ points of $Z\cap Y$. Hence \( (a-1)b\geq ab-1, \) which implies $b\leq 1$. Therefore $b=1$, and $Y$ is a linear space. This concludes the proof.
\end{proof}

\begin{rmk}\label{N34}{\rm 
An immediate application of Lemma \ref{linearX} shows that, for $N=3,4$, every reducible OADP variety whose components are both positive-dimensional is of the form
\[
P\cup Y,
\]
where $P=\p^r\subset\p^N$, with $r=1$ if $N=3$ and $r=1$ or $2$ if $N=4$.

The conclusion of Lemma \ref{linearX} does not extend to the case in which the intersection of the two components has positive dimension. Indeed, let \( Z,Y\subset\p^5 \) be two smooth quadric surfaces such that
\[
\langle Z\rangle\cap\langle Y\rangle
=
Z\cap Y
=
L
\]
is a line. For a general point $q\in\p^5$, \cite[Proposition 9.1.1, third formula]{fulton-intersection} gives
\[
\mult_{\pi_q(L)}
\bigl(\pi_q(Z)\cdot\pi_q(Y)\bigr)
=
3.
\]
It follows that
\[
\alpha_{Z,Y}=4-3=1,
\]
so that $Z\cup Y$ is a reducible OADP variety, although neither component is a linear space.
}
\end{rmk}

We construct several families of reducible OADP surfaces in $\p^5$ of the form
\[
P\cup S,
\]
where $P$ is a plane, generalizing the analogous construction for OADP curves in $\p^3$.

\begin{ex}\label{congP3}
Let $P\subset\p^3$ be a line, and let $Y\subset\p^3$ be a rational curve of degree $d\geq 1$ such that $P\cap Y$ has length $d-1$, or equivalently, $P$ is a $(d-1)$-secant line to $Y$. By Remark \ref{N34}, every reducible OADP curve in $\p^3$ is of this form. If $d\geq 7$, a simple dimension count shows that a general curve $P\cup Y$ of this form is not contained in a cubic surface.

A similar construction applies to surfaces. Let $R\subset\p^{a+1}$ be a smooth surface of minimal degree $a\geq 4$, and let $T=\p^{a-2}$ be the linear span of $a-1$ general points of $R$. Projecting from a general linear subspace \( M=\p^{a-5}\subset T \) gives a smooth surface \( S=\pi_M(R)\subset\p^5 \) and a plane \( P=\pi_M(T)\subset\p^5 \) such that $P\cap S$ consists of the images of the $a-1$ points in $T\cap R$. Thus, $P\cap S$ has length $a-1$, and Remark \ref{intOADP} shows that $P\cup S$ is a reducible OADP surface.

Such a surface can be contained in a cubic fourfold only for $a=4,5$. Indeed, for $a\geq 6$, formula \eqref{defdelta} would give \( \delta=-4a^2+28a-19<21, \) which is impossible. If $a=4$ and $R\subset\p^5$ is a smooth quartic rational normal scroll, then \eqref{defdelta} gives $\delta=29$. If $a=4$ and $R\subset\p^5$ is a Veronese surface, then it gives $\delta=45$. Finally, if $a=5$, it gives $\delta=21$. These constructions recover geometric descriptions of the corresponding well-known divisors $\mathcal S_\delta$.

The divisor $\mathcal S_{21}$ admits another familiar description. If $P,S\subset X$ are two disjoint planes, then \eqref{defdelta} gives $\delta=21$. Thus, $\mathcal S_{21}$ has several different geometric descriptions. A general cubic fourfold $X\in\mathcal S_{21}$ contains different reducible OADP surfaces of the form $P\cup S$, and hence different sections of the associated quadric fibration. Some of these sections are isolated in $X$, whereas others vary in a two-dimensional family, as expected.
\end{ex}

This example shows that the main difficulty is not to construct families of reducible OADP varieties
\[
P\cup S\subset\p^{2r+1},
\qquad
P=\p^r,
\]
of pure dimension $r$, but rather to find such varieties contained in smooth cubic hypersurfaces. The main results of \cite{HassettPlane} ensure their existence for infinitely many values \( \delta\equiv 5\pmod 8, \) and conjecturally for every such $\delta\geq 21$, but no explicit description of the resulting infinite series of rational surfaces is known. In the next section, we construct several new examples.

In Section \ref{SecantLines}, we develop a general framework that illustrates the complexity of the problem and explains why the examples considered here are among the simplest possible ones.

\section{Some reducible OADP surfaces contained in cubic fourfolds}\label{Exseries}

We begin with a family of reducible OADP surfaces parametrized by integers $b\geq 5$.

\begin{ex}\label{Ybspecial}
Let $R_b\subset\p^{b+1}$ be a smooth rational normal scroll of degree $b\geq 5$, let $L\subset R_b$ be a ruling line, and choose $b-3$ general points \( p_1,\ldots,p_{b-3}\in R_b. \) Set
\[
T=\langle L,p_1,\ldots,p_{b-3}\rangle=\p^{b-2},
\]
and let $M=\p^{b-5}\subset T$ be a general linear subspace. Projection from $M$ gives a plane \( P=\pi_M(T)\subset\p^5 \) and a smooth rational scroll \( S_b=\pi_M(R_b)\subset\p^5. \) Writing
\[
L_M=\pi_M(L)
\qquad\text{and}\qquad
q_i=\pi_M(p_i),
\]
we have
\[
S_b\cap P=L_M\cup\{q_1,\ldots,q_{b-3}\}.
\]
We claim that $P\cup S_b\subset\p^5$ is a reducible OADP surface for every $b\geq 5$.

Indeed, let $H_1,H_2\subset\p^5$ be general hyperplanes containing $P$, and write
\[
H_j\cap S_b=L_M+C_j,
\qquad j=1,2.
\]
Then \( C_j\cdot L_M=1 \) and
\[
C_j^2
=
H_j\cdot C_j-L_M\cdot C_j
=
b-1-1
=
b-2.
\]
It follows that \( C_1\cdot C_2=b-2. \) Since the $b-3$ points $q_1,\ldots,q_{b-3}$ belong to $C_1\cap C_2\cap P$, the general linear space
\[
\p^3=H_1\cap H_2
\]
containing $P$ meets $S_b$ in exactly one point outside $P$. Therefore, $S_b$ is a section of the projection from $P$, and Proposition \ref{degpiP} shows that $P\cup S_b$ is a reducible OADP surface.

Suppose first that $S_b$ is contained in a cubic fourfold $X\subset\p^5$. Formula \eqref{selfS} gives \( S_b^2=3b-2, \) and hence
\[
\disc\langle h^2,S_b\rangle
=
3S_b^2-\deg(S_b)^2
=
9b-6-b^2.
\]
For $b\geq 8$, this number is at most $2$, which is impossible for the discriminant of a nonempty Hassett divisor. Direct computations with \emph{Macaulay2} show that $S_b$ is contained in a cubic fourfold for $b=5,6,7$.

Suppose now that
\[
P\cup S_b\subset X\subset\p^5
\]
for some cubic fourfold $X$. We have \( P\cdot S_b=b-3, \) because \( \mult_{L_M}(P\cdot S_b)=0 \) by \cite[Proposition 9.1.1, third formula]{fulton-intersection}; see also \eqref{PdotY}. Formula \eqref{defdelta} then gives \( \delta(S_b)=-4b^2+36b-43. \) The condition \( \delta(S_b)\geq 21 \) forces $b=5$ or $b=6$. For both values, one verifies with \emph{Macaulay2} that $P\cup S_b$ is contained in a cubic fourfold. The corresponding discriminants are
\[
\delta(S_5)=37
\qquad\text{and}\qquad
\delta(S_6)=29.
\]

\end{ex}

We now adapt the preceding construction to produce our main families of reducible OADP surfaces. Some of these surfaces are contained in cubic fourfolds and will be studied in detail below.

\begin{ex}\label{Yer}
Let \( V_e\subset\p^{\frac{e(e+3)}{2}} \) be the Veronese surface defined by forms of degree $e\geq 3$, and let \( 0\leq r<\frac{(e-2)(e+5)}{2}. \) Projecting $V_e$ from the linear span of $r$ general points $p_1,\ldots,p_r\in V_e$, we obtain a smooth surface \( R_{e,r}\subset\p^{\frac{e(e+3)}{2}-r} \) of degree $e^2-r$.

Let $D_{e,r}\subset\p^2$ be a smooth curve of degree $e-1$ passing through the points $p_1,\ldots,p_r$ (such a curve exists if $r\leq\frac{(e-1)(e+2)}{2}$, which we also assume), and denote by $E_{e,r}\subset R_{e,r}$ its image under the Veronese embedding and the subsequent projection. The curve $E_{e,r}$ is smooth, with
\[
\deg(E_{e,r})=e(e-1)-r
\qquad\text{and}\qquad
g(E_{e,r})=\frac{(e-2)(e-3)}{2}.
\]
By the Riemann--Roch theorem,
\[
\dim\langle E_{e,r}\rangle
=
e(e-1)-r-\frac{(e-2)(e-3)}{2}
=
\frac{e^2+3e-6}{2}-r
\geq 3.
\]
Choose a general linear subspace \( M=\p^{\frac{e^2+3e-6}{2}-r-3}\subset\langle E_{e,r}\rangle. \) Projection from $M$ gives a rational map \( \pi_M:\p^{\frac{e(e+3)}{2}-r}\dashrightarrow\p^5. \) Then \( P=\pi_M(\langle E_{e,r}\rangle)\subset\p^5 \) is a plane containing the plane curve
\[
C_{e,r}=\overline{\pi_M(E_{e,r})},
\]
which has degree \( \deg(C_{e,r})=e(e-1)-r \) and geometric genus \( g(C_{e,r})=\frac{(e-2)(e-3)}{2}. \) Since the projection of $E_{e,r}$ from $M$ is general, the curve $C_{e,r}\subset P$ has only ordinary nodes. These are also the nodes of
\[
S_{e,r}=\overline{\pi_M(R_{e,r})}
\]
and exhaust its singular locus. Therefore,
\begin{equation}\label{nodesYer}
\nu(S_{e,r})
=
\frac{\bigl(e(e-1)-r-1\bigr)\bigl(e(e-1)-r-2\bigr)}{2}
-
\frac{(e-2)(e-3)}{2}
>0,
\end{equation}
and
\[
P\cap S_{e,r}=C_{e,r}.
\]

We first claim that $P\cup S_{e,r}\subset\p^5$ is a reducible OADP surface for every $e\geq 2$ and every \( 0\leq r\leq\frac{(e-2)(e+5)}{2}. \) This range includes two classical boundary cases, in which no projection is needed. For $e=2$, one has $r=0$, the surface
\[
S_{2,0}=V_2
\]
is the Veronese surface, and $C_{2,0}$ is a conic spanning the plane $P$. For $e=3$ and $r=4$, the surface $S_{3,4}\subset\p^5$ is a smooth quintic del Pezzo surface, while $C_{3,4}\subset S_{3,4}$ is a conic spanning $P$.

We use the same argument as in Example \ref{Ybspecial}. Let $H_1,H_2\subset\p^5$ be two general hyperplanes containing $P$. Then
\[
H_j\cap S_{e,r}=C_{e,r}+L_j,
\qquad j=1,2,
\]
where $L_j\subset S_{e,r}$ is a rational curve of degree $e$, obtained as the image of a line $l_j\subset\p^2$. Thus,
\[
H_1\cap H_2\cap S_{e,r}\setminus P
=
L_1\cap L_2
\]
consists of the image of the point $l_1\cap l_2$. Hence a general $\p^3$ containing $P$ meets $S_{e,r}$ in exactly one point outside $P$. Therefore, $S_{e,r}$ is a section of the projection from $P$, and Proposition \ref{degpiP} shows that $P\cup S_{e,r}$ is a reducible OADP surface.

We return to the assumptions
\[
e\geq 3
\qquad\text{and}\qquad
0\leq r<\frac{(e-2)(e+5)}{2}
\]
and observe that a necessary condition for $S_{e,r}$ to be contained in a cubic fourfold is that the integer
\begin{equation}\label{bounder}
d_{e,r}
:=
\disc\langle h^2,S_{e,r}\rangle
=
3S_{e,r}^2-\deg(S_{e,r})^2
\end{equation}
satisfies
\[
d_{e,r}>6
\qquad\text{and}\qquad
d_{e,r}\equiv 0,2\pmod 6.
\]
A computation using \eqref{selfS} gives
 \( d_{e,r} = 2e^4-6e^3-4e^2r+9e^2+6er+2r^2-3e-6r+6. \) Most of the surfaces $S_{e,r}$ are not contained in cubic fourfolds, even when they satisfy these necessary numerical conditions. Still fewer reducible OADP surfaces $P\cup S_{e,r}$ are contained in cubic fourfolds. Notice that, if \( \deg(C_{e,r})=e(e-1)-r>3, \) then every cubic hypersurface containing $S_{e,r}$ also contains the plane $P$.

Suppose now that
\[
P\cup S_{e,r}\subset X\subset\p^5
\]
for some cubic fourfold $X$. We compute the discriminant of the lattice \( \langle h^2,P,S_{e,r}\rangle. \) It remains to determine $P\cdot S_{e,r}$.

Let $Q\subset X$ be a general fibre of the quadric fibration induced by projection from $P$. Since $S_{e,r}$ is a section of this fibration, the intersection $Q\cap S_{e,r}$ consists of one point outside $P$, together with the intersection
\[
(Q\cap P)\cap C_{e,r}.
\]
The curve $Q\cap P$ is a conic, while $C_{e,r}\subset P$ has degree $e(e-1)-r$. Therefore, Bézout's theorem gives
\[
Q\cdot S_{e,r}
=
2\bigl(e(e-1)-r\bigr)+1
=
2e^2-2e-2r+1.
\]
Since \( h^2=P+Q, \) it follows that
\[
P\cdot S_{e,r}
=
h^2\cdot S_{e,r}-Q\cdot S_{e,r}
=
e^2-r-2e(e-1)+2r-1
=
r+2e-e^2-1.
\]
Consequently,
\begin{equation}\label{deltaer}
\begin{split}
\delta(S_{e,r})
&=
\det
\begin{pmatrix}
3 & 1 & e^2-r\\
1 & 3 & r+2e-e^2-1\\
e^2-r & r+2e-e^2-1 & S_{e,r}^2
\end{pmatrix}\\
&=
4e^2+4e-8r+13.
\end{split}
\end{equation}

More generally, suppose that
\[
P\cup S\subset X\subset\p^5
\]
is a reducible OADP surface contained in a cubic fourfold, and that the scheme-theoretic intersection $P\cap S$ consists of a curve $C$ and a zero-dimensional residual scheme. Since $S$ is a section of the quadric fibration induced by projection from $P$, we have \( Q\cdot S=2\deg(C)+1. \) Hence
\begin{equation}\label{PdotY}
P\cdot S
=
h^2\cdot S-Q\cdot S
=
\deg(S)-2\deg(C)-1.
\end{equation}

The construction above admits the following generalization. Let
\[
R\simeq\Bl_{p_1,\ldots,p_r}(\p^2)\subset\p^s,
\qquad s\geq 6,
\]
be a smooth rational surface embedded by the complete linear system
\[
|D|
=
\left|eL-\sum_{i=1}^r m_iE_i\right|,
\]
where $e\geq 3$ and \( 1\leq m_1\leq\cdots\leq m_r. \) Assume that $R$ contains an irreducible curve $E$ of degree $c\geq 3$ such that \( \langle E\rangle=\p^{s-3}, \) and choose a general linear subspace \( M=\p^{s-6}\subset\langle E\rangle. \) Suppose that $|D-E|$ defines a birational map onto $\p^2$, or equivalently that the corresponding plane linear system is a homaloidal net. Then, arguing as above,
\[
S=\pi_M(R)
\qquad\text{and}\qquad
P=\pi_M(\langle E\rangle)
\]
form a reducible OADP surface $P\cup S\subset\p^5$. The surface $S$ has \( \nu = \frac{(c-1)(c-2)}{2}-g(E) \) nodes, and its intersection with $P$ consists of the curve
\[
C=\pi_M(E),
\qquad
\deg(C)=c,
\]
possibly together with a zero-dimensional residual scheme.

If $P\cup S$ is contained in a cubic fourfold, then
\begin{equation}\label{deltaex}
\delta
=
\det
\begin{pmatrix}
3 & 1 & \deg(S)\\
1 & 3 & \deg(S)-2c-1\\
\deg(S) & \deg(S)-2c-1 & S^2
\end{pmatrix}.
\end{equation}
\end{ex}
\medskip

Applying the construction above, we obtained the examples summarized in Table \ref{TABex}. We encode the pair of plane linear systems defining $R$ and $E$ by
\begin{equation}\label{notationSC}
\mathfrak{S}\left(
\begin{array}{ccccc}
e & a_1 & a_2 & a_3 & \cdots\\
\alpha & b_1 & b_2 & b_3 & \cdots
\end{array}
\right),
\end{equation}
where $e$ and $\alpha$ are their degrees, $a_i$ is the number of base points of multiplicity $i$ of the first system, and $b_i\leq a_i$ is the number of these points that are simple base points of the second system. Thus,
\[
r=\sum_{i\geq1}a_i
\qquad\text{and}\qquad
c=\deg(C)=e\alpha-\sum_{i\geq1}i\,b_i,
\]
and trailing zero entries may be omitted.

For instance, Example 17 in Table \ref{TABex} has type
\[
\mathfrak{S}\left(
\begin{array}{ccc}
5 & 6 & 2\\
3 & 5 & 2
\end{array}
\right).
\]
Hence the first system has degree $5$, with six simple and two double base points, whereas the second consists of plane cubics passing simply through five of the simple base points and through the two double base points of the first system. In particular, $r=8$. The resulting surface $S\subset\p^5$ has degree $11$ and sectional genus $4$, while $C=P\cap S$ has degree $6$ and $9$ nodes.

\subsection{The inverse of the projection from $P$}\label{inverseProj}

Let $S\subset\p^5$ be one of the surfaces constructed in Example \ref{Yer}, and let \( f:\Bl_{p_1,\ldots,p_r}(\p^2)\longrightarrow S\subset\p^5 \) be the finite birational morphism induced by the projection from $M$. Then $\Sing(S)=\Sing(C)\subset P$, where $C=P\cap S$, and each node of $S$ has exactly two preimages under $f$.

The hyperplanes containing $P$ cut out on $S$ the curve $C$ together with the curves of the homaloidal net $|D-E|$. Hence the composition of $f$ with the projection from $P$ is the birational morphism \( \Bl_{p_1,\ldots,p_r}(\p^2)\to\p^2 \) defined by $|D-E|$, that is, the composition of the blow-up morphism with the Cremona transformation \( \phi:\p^2\dashrightarrow\p^2 \) defined by the corresponding plane linear system. It follows that the inverse \( s:\p^2\dashrightarrow S \) of the birational map $\pi_P|_S:S\dashrightarrow\p^2$ is undefined exactly at the images $q_1,\ldots,q_r\in\p^2$ of the curves contracted by this morphism, and that $E$ is mapped onto the plane curve \( \overline C=\overline{\phi(E)}\subset\p^2. \) Moreover, $s$ maps $\overline C$ birationally onto $C$, and the two points of $\overline C$ corresponding to the two preimages of a node of $S$ are both mapped to that node.

For the surfaces $S_{e,r}$, one has $D-E=L$. Thus, $\phi$ is the identity, $q_i=p_i$, the curve $\overline C$ is the curve $D_{e,r}$, and $s$ is defined by a linear subsystem of the system of plane curves of degree $e$ passing through $p_1,\ldots,p_r$, that is, by a linear system of the same type as the one defining $f$. For the other surfaces of Example \ref{Yer}, $\phi$ is a nontrivial Cremona transformation, and the types of $\overline C$ and of the linear system defining $s$ are obtained from those of $E$ and $|D|$ through the action of $\phi$.

\section{Parameter count and construction of some divisors $\mathcal S_\delta\subset \mathcal C_8$}\label{ParCount}

We return to the examples in Table \ref{TABex}. For a type $\mathfrak S$ as in \eqref{notationSC}, let $p_1,\ldots,p_r\in\p^2$ be general points, and denote by $m_j$ and $n_j$ the multiplicities of $p_j$ in the first and second plane linear systems, respectively. Consider the blow-up $\Bl_{p_1,\ldots,p_r}(\p^2)$, embedded in $\p^s$ by the very ample complete linear system
\[
\left|eL-\sum_{j=1}^r m_jE_j\right|,
\]
where $L$ is the pullback of the class of a line in $\p^2$ and $E_j$ is the exceptional divisor over $p_j$. Let
\[
E\in\left|\alpha L-\sum_{j=1}^r n_jE_j\right|
\]
be a smooth irreducible curve whose linear span has dimension $s-3$, and choose a general linear space \( M=\p^{s-6}\subset\langle E\rangle. \) Projection from $M$ gives
\[
S=\pi_M\bigl(\Bl_{p_1,\ldots,p_r}(\p^2)\bigr)\subset\p^5,
\qquad
P=\pi_M(\langle E\rangle),
\qquad
C=\pi_M(E)\subset S.
\]
In all the examples considered here, $P\cap S=C$, and $S$ has \( \frac{(c-1)(c-2)}{2}-g(E) \) nodes, where $c=\deg(C)=e\alpha-\sum_{j=1}^r m_jn_j$; these nodes are precisely the nodes of $C$.

We say that the resulting pair $(S,C)$, or equivalently $(S,P)$, is of type $\mathfrak S$. We now count the dimension of the family of such pairs.

\begin{Proposition}\label{moduliSC}
The pairs $(S,C)$ of type $\mathfrak S$ depend on
\[
n(\mathfrak S):=
\max\{0,2(r-4)\}
+\dim\left|\alpha L-\sum_{j=1}^r n_jE_j\right|
+3(s-5)+35-\max\{0,8-2r\}
\]
parameters.
\end{Proposition}

\begin{proof}
Up to projective transformations of $\p^2$, the points $p_1,\ldots,p_r$ depend on
\[
\max\{0,2(r-4)\}
\]
parameters, while the curve $E$ varies in a linear system of dimension
\[
\dim\left|\alpha L-\sum_{j=1}^r n_jE_j\right|.
\]
The centers $M=\p^{s-6}\subset\langle E\rangle\simeq\p^{s-3}$ form a Grassmannian of dimension \( \dim\mathbb G(s-6,s-3)=3(s-5). \) Finally, we add the $35$ parameters corresponding to $\PGL{6}$ and subtract the $\max\{0,8-2r\}$ parameters corresponding to the automorphisms of $\Bl_{p_1,\ldots,p_r}(\p^2)$. This gives the stated formula.
\end{proof}

Let $\mathcal H_{\mathfrak S}$ be the closure of the family of pairs $(S,P)$ of type $\mathfrak S$. It is irreducible and unirational, being dominated by the rational parameter space used in the proof of Proposition \ref{moduliSC}, and it has dimension $n(\mathfrak S)$. Denoting by $\widetilde{\mathcal C}$ the parameter space of cubic hypersurfaces in $\p^5$, consider the incidence variety
\[
\mathcal I_{\mathfrak S}
=
\left\{((S,P),X)\in
\mathcal H_{\mathfrak S}\times\widetilde{\mathcal C}
:\ S\cup P\subset X\right\}.
\]
By semicontinuity, $\h^0(\mathcal I_{S\cup P}(3))$ attains its minimum value on a dense open subset of $\mathcal H_{\mathfrak S}$, over which $\mathcal I_{\mathfrak S}$ is a projective bundle. In what follows, we replace $\mathcal I_{\mathfrak S}$ by the closure of this projective bundle, namely the irreducible component of $\mathcal I_{\mathfrak S}$ dominating $\mathcal H_{\mathfrak S}$. It is unirational, since $\mathcal H_{\mathfrak S}$ is, and for a general pair $(S,P)\in\mathcal H_{\mathfrak S}$ its dimension is
\[
\dim\mathcal I_{\mathfrak S}
=
n(\mathfrak S)+\h^0(\mathcal I_{S\cup P}(3))-1.
\]

Since a plane $P\subset X$ does not deform inside $X$, the fibre of the projection $\mathcal I_{\mathfrak S}\to\widetilde{\mathcal C}$ through a general point $((S,P),X)$ has dimension at most $\h^0(N_{S/X})$, the dimension of the tangent space at $[S]$ to the Hilbert scheme of $X$. Assume that a general cubic fourfold containing a general pair $S\cup P$ is smooth, and let $\mathcal C_{\mathfrak S}\subseteq\mathcal C_8$ be the closure of the image of $\mathcal I_{\mathfrak S}$, that is, of the locus of cubic fourfolds containing $S\cup P$ for some $(S,P)\in\mathcal H_{\mathfrak S}$. It is irreducible and unirational, since $\mathcal I_{\mathfrak S}$ is. The locus $\widetilde{\mathcal C}_8\subset\widetilde{\mathcal C}$ of cubic fourfolds containing a plane has dimension $54$. Hence, by the fibre dimension theorem and after quotienting by $\PGL{6}$,
\begin{equation}\label{eq:stima}
\codim_{\mathcal C_8}\mathcal C_{\mathfrak S}
\leq
54-
\left(
n(\mathfrak S)
+\h^0(\mathcal I_{S\cup P}(3))
-1
-\h^0(N_{S/X})
\right)
\end{equation}
for a general $(S,P)\in\mathcal H_{\mathfrak S}$ and a general cubic fourfold $X$ containing $S\cup P$.

This gives a practical criterion for determining whether cubic fourfolds containing a general pair $(S,P)$ of type $\mathfrak S$ fill a divisor $\mathcal S_\delta$.

\begin{Proposition}\label{suffSd}
Let $(S,P)\in\mathcal H_{\mathfrak S}$ be a general pair, and let $X$ be a general cubic fourfold containing $S\cup P$. Assume that $X$ is smooth, that $\delta=\disc\langle h^2,P,S\rangle$, and that
\[
n(\mathfrak S)
+\h^0(\mathcal I_{S\cup P}(3))
-\h^0(N_{S/X})
=54.
\]
Then $\mathcal C_{\mathfrak S}=\mathcal S_\delta$. In particular, $\mathcal S_\delta$ is unirational.
\end{Proposition}

\begin{proof}
By the inequality above, $\mathcal C_{\mathfrak S}$ is an irreducible subvariety of codimension at most one in $\mathcal C_8$. A general member $[X']$ of $\mathcal C_{\mathfrak S}$ is a smooth cubic fourfold containing $S'\cup P'$ for a general pair $(S',P')\in\mathcal H_{\mathfrak S}$. Thus, $S'$ is a section of the quadric fibration of $X'$ determined by $P'$, by Proposition \ref{degpiP}, and the lattice $\langle h^2,P',S'\rangle\subseteq\H^{2,2}(X',\mathbb Z)$ has discriminant $\delta$, since the relevant intersection numbers are constant in the family. If $Q'$ is a fibre of this quadric fibration, then $\langle h^2,P',S'\rangle=\langle h^2,Q',S'\rangle$, and $Q'\cdot S'$ is odd by Remark \ref{TvsY}. Hence $[X']\in\mathcal S_\delta$. Since $\mathcal S_\delta$ is an irreducible divisor in $\mathcal C_8$, it follows that $\mathcal C_{\mathfrak S}=\mathcal S_\delta$, which is therefore unirational.
\end{proof}

\begin{rmk}\label{explicitValues}{\rm
In practice, the integers $\h^0(\mathcal I_{S\cup P}(3))$ and $\h^0(N_{S/X})$ are computed for an explicit pair $(S_0,P_0)$ and an explicit cubic fourfold $X_0$. By semicontinuity, $\h^0(N_{S_0/X_0})$ is an upper bound for the value of $\h^0(N_{S/X})$ at a general point, so that it can be used in the inequality above. Similarly, $\h^0(\mathcal I_{S_0\cup P_0}(3))$ is an upper bound for the general value of $\h^0(\mathcal I_{S\cup P}(3))$. On the other hand, the latter is at least $56-\h^0(\mathcal O_{S\cup P}(3))$ for a general pair, and hence at least $56-\h^0(\mathcal O_{S_0\cup P_0}(3))$. Therefore, the two values coincide whenever $S_0\cup P_0$ imposes independent conditions on cubic hypersurfaces.
}
\end{rmk}

\begin{rmk}\label{refinedCount}{\rm
The bound $\h^0(N_{S/X})$ for the dimension of the fibre of $\mathcal I_{\mathfrak S}\to\widetilde{\mathcal C}$ can be improved. Since $P$ does not deform inside $X$, near a general point $((S,P),X)$ this fibre consists of pairs $(S',P)$ with the same plane $P$, and $S'$ is the closure of $(S'\cup P)\setminus P$. Hence the fibre embeds into the closed subscheme of the Hilbert scheme of $X$ parametrizing the subschemes that contain $P$, whose tangent space at $[S\cup P]$ is
\[
\operatorname{Hom}(\mathcal I_{S\cup P/X},\mathcal I_{P/S\cup P}),
\]
where $\mathcal I_{S\cup P/X}$ is the ideal sheaf of $S\cup P$ in $X$ and $\mathcal I_{P/S\cup P}$ is that of $P$ in $S\cup P$. This is the space of first-order deformations of $S\cup P$ inside $X$ that still contain $P$; note that the curve $S\cap P$ is not required to stay fixed. Denoting by $t$ its dimension, $\h^0(N_{S/X})$ can be replaced by $t$ in the inequality \eqref{eq:stima} and in Proposition \ref{suffSd}; as in Remark \ref{explicitValues}, the value of $t$ computed for an explicit pair is an upper bound for its general value. For all the examples in Table \ref{TABex}, one has $t=\h^0(N_{S/X})$. For Examples 1--12 of Table \ref{TABex2}, instead, one has $t<\h^0(N_{S/X})$, and the resulting upper bound for $\codim_{\mathcal C_8}\mathcal C_{\mathfrak S}$ is $1$, so that the corresponding cubic fourfolds fill the divisors $\mathcal S_\delta$. In particular, Example 12 of Table \ref{TABex2} shows that $\mathcal S_{77}$ is unirational. The tables in Section \ref{SummaryTables} were computed using $\h^0(N_{S/X})$; see also Section \ref{Computations}.
}
\end{rmk}

We illustrate the method with a relevant example; see Section \ref{Computations} for the implementation of this parameter count in the package \emph{SpecialFanoFourfolds}.

\subsection{The divisor $\mathcal S_{53}$}\label{S53}

Let $S\subset\p^5$ be a surface of type
\[
\mathfrak{S}\left(\begin{array}{ccccc}
4 & 5 & 0 & 0 & 0 \\
3 & 5 & 0 & 0 & 0
\end{array}\right),
\]
that is, the surface $S_{4,5}$ of Example \ref{Yer}; see also Example 18 in Table \ref{TABex}. In this case, $C=S\cap P$ is a curve of degree $7$ and geometric genus $1$ with $14$ nodes, and these nodes constitute the entire singular locus of $S$. A general cubic fourfold $X\subset\p^5$ containing $S\cup P$ is smooth, formula \eqref{deltaex} gives $\delta=53$, and formula \eqref{selfS} shows that $[X]\in\mathcal C_{86}$.

Here $r=5$ and $s=9$, so that
\[
n(\mathfrak S)=2+4+12+35=53.
\]
Direct computations in \emph{Macaulay2} yield
\[
\h^0(\mathcal I_{S\cup P}(3))=6
\qquad\text{and}\qquad
\h^0(N_{S/X})=5.
\]
Therefore,
\[
n(\mathfrak S)
+\h^0(\mathcal I_{S\cup P}(3))
-\h^0(N_{S/X})
=
53+6-5=54.
\]
It follows from Proposition \ref{suffSd} that the cubic fourfolds containing pairs $(S,P)$ of the above type fill $\mathcal S_{53}$, which is therefore unirational. By Section \ref{SdeltaSec}, $\mathcal S_{53}$ is an irreducible component of $\mathcal C_8\cap\mathcal C_{86}$; thus, it is a unirational codimension-two locus of rational cubic fourfolds contained in the admissible divisor $\mathcal C_{86}$, which is not uniruled by \cite{TaVaA}. Examples 13--17 of Table \ref{TABex} give other geometric descriptions of $\mathcal S_{53}$, in which the class of $S$ has discriminant $48$ or $56$.

\section{Secant lines to $P\cup S$ and the associated K3 surface}\label{SecantLines}

A cubic fourfold containing a plane carries a natural quadric fibration, and the lines contained in its fibres are parametrized by a relative Hilbert scheme, which is a $\p^1$-bundle over the associated K3 surface. If the cubic fourfold contains a reducible OADP surface $P\cup S$, the secant lines to $P\cup S$ contained in it form a surface $U_S$ in this relative Hilbert scheme, which is birational to the K3 surface. In this section, we describe $U_S$ and show how the Fano maps of $P\cup S$ produce projective models of the K3 surface.

\subsection{The relative Hilbert scheme of lines and the surface $U_S$}\label{associatedK3}

Let $X\subset\p^5$ be a cubic fourfold containing a plane $P$, and let \( f_P:\Bl_P(X)\longrightarrow\p^2 \) be the quadric fibration induced by projection from $P$.

We say that the pair $(X,P)$ is \emph{general} if the discriminant sextic $\Delta_X\subset\p^2$ of $f_P$ is smooth. A general cubic fourfold $X\in\mathcal C_8$ contains a unique plane $P$, and the cubic fourfolds $X$ for which $(X,P)$ is general form a nonempty open subset of $\mathcal C_8$.

If $X\in\mathcal S_\delta$, we say that $X$ is general if it contains a plane $P$ such that $(X,P)$ is general. If this condition holds for one element of $\mathcal S_\delta$, then it holds for a general element.

For a general pair $(X,P)$, let $K_X$ be the double cover of $\p^2$ branched along $\Delta_X$. Then $K_X$ is a smooth K3 surface endowed with a polarization $D_X$ of degree $2$, and the linear system $|D_X|$ defines the double cover \( \alpha_X:K_X\longrightarrow\p^2. \) We call the polarized pair $(K_X,D_X)$ the \emph{K3 surface associated with $(X,P)$}.

If $X\in\mathcal C_8$ is Noether--Lefschetz general, that is, if $\rk\H^{2,2}(X,\mathbb Z)=2$, then $(K_X,D_X)$ is very general in the moduli space of polarized K3 surfaces of degree $2$. Equivalently,
\[
\Pic(K_X)=\langle D_X\rangle\simeq\mathbb Z D_X.
\]

We use the following well-known fact; see, for example, \cite{HassettPlane,kuz4fold,LevicoK}.

\begin{Proposition}\label{HilbPX}
Let $(X,P)$ be a general pair. Then the relative Hilbert scheme $\mathcal H_{P,X}$ of lines contained in the fibres of $f_P$ is a $\p^1$-bundle over $K_X$.
\end{Proposition}

\begin{proof}
For a point $z\in\p^2$, the Hilbert scheme of lines contained in the quadric surface $f_P^{-1}(z)$ is the disjoint union of two copies of $\p^1$ if $z\notin\Delta_X$, corresponding to the two rulings of a smooth quadric surface, and is a single copy of $\p^1$ if $z\in\Delta_X$. The two rulings are separated by the double cover $K_X\to\p^2$. Consequently, $\mathcal H_{P,X}$ is a smooth irreducible threefold endowed with a morphism \( \beta_X:\mathcal H_{P,X}\longrightarrow K_X \) whose fibres are isomorphic to $\p^1$.
\end{proof}

\medskip

Let $S\subset X$ be a section of $\pi_P|_X:X\dashrightarrow\p^2$. We denote by $U_S\subset\mathcal H_{P,X}$ the closure of the locus parametrizing the \emph{proper secant lines} to $P\cup S$ contained in $X$, namely, the lines in $X$ joining a general point $y\in S$ to a point $p\in P$.

\begin{cor}\label{Usm}
Let $(X,P)$ be a general pair, and let $S\subset X$ be a section of $\pi_P|_X$. Then $U_S$ is an irreducible surface, and \( \beta_X|_{U_S}:U_S\longrightarrow K_X \) is a birational morphism whose positive-dimensional fibres are fibres of $\beta_X$.
\end{cor}

\begin{proof}
Let $y\in S$ be a general point and set $z=\pi_P(y)$. A proper secant line $\langle p,y\rangle$ to $P\cup S$ contained in $X$ projects from $P$ to $z$ and is therefore contained in the smooth quadric fibre $Q_z=f_P^{-1}(z)$. It determines a point of $U_S$ lying over the point of $\alpha_X^{-1}(z)$ corresponding to the ruling of $Q_z$ containing $\langle p,y\rangle$.

Through $y$ there pass exactly two lines contained in $Q_z$ and meeting $P$, one in each ruling. By the generality of $y$, neither line is contained in $S$. Thus, each of the two points of $\alpha_X^{-1}(z)$ has a unique preimage in $U_S$, and $\beta_X|_{U_S}$ is birational. The surface $U_S$ is irreducible by construction.

Let $s\in K_X$ be one of the finitely many points at which the inverse of $\beta_X|_{U_S}$ is not defined. By Zariski's Main Theorem, the fibre of $\beta_X|_{U_S}$ over $s$ is connected and one-dimensional. Since it is contained in $\beta_X^{-1}(s)\simeq\p^1$, it coincides with $\beta_X^{-1}(s)$.
\end{proof}

In general, $U_S$ is neither smooth nor normal. Its structure depends on the geometry of $S$, and it is described in Proposition \ref{doublecoverT} below for the sections considered in this paper.

Let $\mathcal L_{P,X}\to\mathcal H_{P,X}$ be the universal family of lines. The natural morphism $\mathcal L_{P,X}\to\Bl_P(X)$ is generically finite of degree two, since exactly two lines contained in the corresponding quadric fibre pass through a general point of $\Bl_P(X)$. Let $E_S\subset\Bl_P(X)$ be the image of the restriction of $\mathcal L_{P,X}$ to $U_S$, that is, the divisor swept out by the lines parametrized by $U_S$. Then $E_S$ is irreducible and singular along $S$, and the restriction of $\mathcal L_{P,X}$ to $U_S$ maps birationally onto $E_S$. Indeed, the general fibre of the dominant morphism $f_P|_{E_S}:E_S\to\p^2$ is one-dimensional and consists of the two secant lines to $P\cup S$ passing through the corresponding point of $S$.

By a slight abuse of notation, we shall also denote by $E_S$ its image in $X$ and its strict transform in $\Bl_{P\cup S}(X)$.

\subsection{Simple sections}

We work under the following assumption: the section $S\subset X\subset\p^5$ of $\pi_P|_X$ has at most nodes as singularities, all contained in $P$, and these coincide with the nodes of the plane curve $C=P\cap S$.

When this assumption is satisfied, we shall simply say that $S$ \emph{has at most nodes as singularities}. This is a natural hypothesis for several reasons, as we shall see in the following sections. Moreover, most of the surfaces used in the literature to describe Hassett divisors $\mathcal C_d\subset\mathcal C$ geometrically have singularities of this type.

Let \( \pi:\Bl_P(X)\longrightarrow X \) be the blow-up morphism. Under the above assumption, the strict transform $T\subset\Bl_P(X)$ of $S$ is smooth, since blowing up $P$ resolves the singularities of $S$. Moreover, $T$
is endowed with two birational morphisms:
\[
f_T:=f_P|_T:T\longrightarrow\p^2
\qquad\text{and}\qquad
\pi|_T:T\longrightarrow S.
\]
The latter is finite and birational, so that $T$ is naturally identified with the smooth normalization of $S$.

Let $p_1,\ldots,p_r\in\p^2$, with $r\geq 0$, be the points at which the rational section \( s:\p^2\dashrightarrow S\subset\p^5 \) is not defined; they are also the proper base points of the inverse rational map \( f_T^{-1}:\p^2\dashrightarrow T. \) When $S$ is one of the surfaces of Example \ref{Yer}, these are the points denoted by $q_1,\ldots,q_r$ in Section \ref{inverseProj}.

\begin{defin}\label{simpleSection}
A section $S\subset X\subset\p^5$ having at most nodes as singularities is called \emph{simple} if the birational morphism \( f_T:T\longrightarrow\p^2 \) is isomorphic to the blow-up \( \Bl_{p_1,\ldots,p_r}(\p^2)\longrightarrow\p^2. \) \end{defin}

\begin{rmk}\label{simplevsnonsimple}
{\rm 
Since $f_T:T\to\p^2$ is a birational morphism between smooth surfaces, it is a composition of blow-ups at points. Thus, each fibre $f_T^{-1}(p_i)$ is a tree of transversely intersecting rational curves contained in the quadric \( Q_i=f_P^{-1}(p_i). \) The section $S$ is simple if and only if
\[
f_T^{-1}(p_i)\simeq\p^1
\]
for every $i=1,\ldots,r$, or equivalently, if the base locus of $f_T^{-1}$ has no infinitely near points.

Among the sections considered here, simple sections are both the easiest to describe and the most common. The constructions of Section \ref{Exseries} provide many examples of simple sections; in fact, all the examples listed in Table \ref{TABex} are simple.

If $r=0$, that is, if $f_T:T\to\p^2$ is an isomorphism, then $U_S\simeq K_X$; see Proposition \ref{doublecoverT} below. This is the case in Examples 3 and 22 of Table \ref{TABex}. In Example 3, the surface $S$ is a Veronese surface and $P$ is the plane spanned by a conic on $S$.
}
\end{rmk}

All the sections considered in our examples are simple, and their base points do not lie on the discriminant curve $\Delta_X$. Under these assumptions, the surface $U_S$ can be described as follows.

\begin{Proposition}\label{doublecoverT}
Let $(X,P)$ be a general pair, and let $S\subset X$ be a simple section whose base points $p_1,\ldots,p_r$ do not lie on $\Delta_X$. For each $i$, the curve $C_i=f_T^{-1}(p_i)$ has bidegree $(1,e_i)$, for some $e_i\geq 0$, on the smooth quadric surface $Q_i=f_P^{-1}(p_i)\simeq\p^1\times\p^1$. Write \( \alpha_X^{-1}(p_i)=\{s_i,s'_i\}, \) where $s_i$ corresponds to a ruling of $Q_i$ whose lines meet $C_i$ in one point. Then the following hold.
\begin{enumerate}
\item[{\rm (i)}] The normalization of $U_S$ is the blow-up $\overline U_S$ of $K_X$ at the points $s_1,\ldots,s_r$ and at the points $s'_i$ with $e_i\geq 1$; more precisely, the composition of the normalization morphism $\overline U_S\to U_S$ with $\beta_X|_{U_S}$ is the blow-up morphism.

\item[{\rm (ii)}] The surface $U_S$ is smooth if and only if $e_i\leq 1$ for every $i$. In this case, $U_S=\overline U_S$ is the blow-up of $K_X$ at $r+j$ points, where $j$ is the number of indices $i$ with $e_i=1$. If $e_i\geq 2$, then $U_S$ is not normal along $\beta_X^{-1}(s'_i)$.
\end{enumerate}
\end{Proposition}

\begin{proof}
Since $C_i\simeq\p^1$ is a smooth curve on $Q_i\simeq\p^1\times\p^1$, its bidegree is $(1,e_i)$ up to exchanging the two factors, and the lines of the two rulings meet $C_i$ in $1$ and $e_i$ points, respectively.

Let $V_S\subset\mathcal H_{P,X}\times T$ be the surface of pairs $([\ell],t)$ such that $t\in\ell$, that is, the preimage of $T$ in the universal family $\mathcal L_{P,X}$. We first show that $V_S$ is the blow-up of $K_X$ at the $2r$ points $s_1,s'_1,\ldots,s_r,s'_r$. Let $B\subset\Bl_P(X)$ be the complement of the vertices of the singular fibres of $f_P$. Then $f_P|_B$ is smooth, and hence so is the projection $K_X\times_{\p^2}B\to K_X$; in particular, $K_X\times_{\p^2}B$ is smooth. The morphism from the preimage of $B$ in $\mathcal L_{P,X}$ to $K_X\times_{\p^2}B$ given by \( ([\ell],x)\mapsto(\beta_X([\ell]),x) \) is bijective, because through a point $x\in B$ there passes exactly one line of each ruling of the fibre containing $x$. Hence it is an isomorphism by Zariski's Main Theorem. Moreover, $T\subset B$: if $t\in T$ lies over a point of $\Delta_X$, then $f_T$ is a local isomorphism at $t$, because $p_i\notin\Delta_X$, and therefore $f_P$ is smooth at $t$. It follows that \( V_S\simeq K_X\times_{\p^2}T. \) Since $f_T$ is the blow-up of $\p^2$ at $p_1,\ldots,p_r$ and $\alpha_X$ is étale over these points, $K_X\times_{\p^2}T$ is the blow-up of $K_X$ at $s_1,s'_1,\ldots,s_r,s'_r$. Denote by $E_i$ and $E'_i$ the exceptional curves over $s_i$ and $s'_i$, respectively. The curve $E_i$ (resp.\ $E'_i$) consists of the pairs $([\ell],t)$ such that $t\in C_i$ and $\ell$ is the line through $t$ in the ruling corresponding to $s_i$ (resp.\ $s'_i$).

Let $\eta:V_S\to\mathcal H_{P,X}$ be the projection. Its image is irreducible and contains the proper secant lines to $P\cup S$ contained in $X$, so that $\eta(V_S)=U_S$. If $f_T(t)\notin\{p_1,\ldots,p_r\}$, the fibre of $f_P$ containing $t$ meets $T$ only at $t$, and so does every line through $t$ in this fibre. Since $\beta_X\circ\eta$ is the blow-up morphism $V_S\to K_X$, it follows that $\eta$ is birational, and that it is injective and immersive away from the exceptional curves. Moreover, $\eta$ maps $E_i$ isomorphically onto $\beta_X^{-1}(s_i)$, and it maps $E'_i$ isomorphically onto $\beta_X^{-1}(s'_i)$ if $e_i=1$, onto $\beta_X^{-1}(s'_i)$ with degree $e_i$ if $e_i\geq 2$, and onto the point $[C_i]$ if $e_i=0$. If $\eta$ maps an exceptional curve isomorphically onto a fibre of $\beta_X$, then $\eta$ is immersive along it, because the kernel of the differential of the blow-up morphism at a point of an exceptional curve is the tangent line to that curve. Consequently, $\eta$ restricts to an isomorphism from \( V_S\setminus\bigcup_{e_i\neq 1}E'_i \) onto its image, which is open in $U_S$, since it is a proper, injective and immersive morphism onto this image.

If $e_i=0$, the fibre of $\beta_X|_{U_S}$ over $s'_i$ is the point $[C_i]$, so that $\beta_X|_{U_S}$ is an isomorphism over a neighbourhood of $s'_i$ by Zariski's Main Theorem. Hence $\eta$ factors through the contraction of the curves $E'_i$ with $e_i=0$, which is the blow-up $\overline U_S$ described in (i). The induced morphism $\overline U_S\to U_S$ is finite and birational, and it is an isomorphism away from the curves $E'_i$ with $e_i\geq 2$. Since $\overline U_S$ is smooth, it is the normalization of $U_S$, which proves (i). If $e_i\geq 2$, the normalization morphism maps $E'_i$ onto $\beta_X^{-1}(s'_i)$ with degree $e_i$, so that it is not an isomorphism over this curve, and $U_S$ is not normal along it. Together with (i), this proves (ii).
\end{proof}

We also denote by $E_i$ and $E'_i$ the exceptional curves of the blow-up $\overline U_S\to K_X$ over $s_i$ and $s'_i$, respectively, the latter only when $e_i\geq 1$. By the proof above, the normalization morphism is an isomorphism over the complement of the curves $\beta_X^{-1}(s'_i)$ with $e_i\geq 2$, and hence $U_S$ is smooth there. In particular, the normalization morphism maps $E_i$ isomorphically onto $\beta_X^{-1}(s_i)$, and $E'_i$ isomorphically onto $\beta_X^{-1}(s'_i)$ when $e_i=1$; in these cases, we identify these curves with their images in $U_S$. Moreover, the lines parametrized by $E_i$ form one of the two rulings of $Q_i$, so that $Q_i\subset E_S$ for every $i$.

For the surfaces of Example \ref{Yer}, the curves $\pi(C_i)\subset S$ are the curves contracted by the projection from $P$, that is, the images in $S$ of the curves contracted by the morphism $\Bl_{p_1,\ldots,p_r}(\p^2)\to\p^2$ defined by $|D-E|$; see Section \ref{inverseProj}. Since $\pi(C_i)$ has bidegree $(1,e_i)$ on the quadric surface $Q_i\subset\langle Q_i\rangle\simeq\p^3$, its degree is $e_i+1$. This allows one to compute the integers $e_i$ from the type $\mathfrak S$. In Examples 1, 2, 4, 6, 9, 12, 15 and 18 of Table \ref{TABex}, one has $e_i=0$ for every $i$, while $r=0$ in Examples 3 and 22. Moreover, $j=1$ in Examples 7, 13 and 21, $j=2$ in Examples 5, 8, 11, 14 and 17, $j=3$ in Examples 20 and 24, and $j=5$ in Examples 19 and 23. Finally, in Examples 8, 10, 14 and 16 there is exactly one index $i$ with $e_i=2$, and $e_i\leq 1$ for all the other indices; by Proposition \ref{doublecoverT}, the surface $U_S$ is not normal in these four examples.

\subsection{Fano maps and projective models of $K_X$}\label{Fanomapmu}

In this subsection, we do not assume that $S\subset X$ is simple. Since a general secant line to $P\cup S$ is not contained in $X$, sending a general pair $(p,y)\in P\times S$ to the residual point of intersection of the line $\langle p,y\rangle$ with $X$ defines a birational map $P\times S\dashrightarrow X$. Its inverse contracts the divisor $E_S\subset X$ onto the image of the birational map $U_S\dashrightarrow P\times S$ sending $[\langle p,y\rangle]$ to $(p,y)$, which is in general a highly singular surface; the general fibres of this contraction are the secant lines to $P\cup S$ contained in $X$.

Since $P\cup S\subset\p^5$ is a reducible OADP surface, through a general point $q\in\p^5$ there passes a unique secant line $\langle p,y\rangle$ to $P\cup S$. Sending $q$ to $(p,y)$ defines a rational map $\p^5\dashrightarrow P\times S$ whose general fibres are the secant lines to $P\cup S$. Thus, the OADP property realizes the abstract join of $P$ and $S$ geometrically, and the restriction of this map to $X$ is the inverse of the birational map $P\times S\dashrightarrow X$ above.

To realize this map by a rational map to a projective space, one considers linear systems of forms of degree $d\geq 2$ on $\p^5$, or more naturally on $\Bl_{P\cup S}(\p^5)$, with multiplicity at least $a\geq 1$ along $P$ and at least $b\geq 1$ along $S$, where $a+b=d$. We also require the corresponding projective linear system to have dimension at least $4$. Finding such linear systems is generally nontrivial. Since $S\subset\p^5$ is nondegenerate in most of our examples and we are interested in surfaces contained in cubic fourfolds, the most favorable case is $a=d-1$ and $b=1$. When $b>1$, the conditions imposed along $S$ are generally too numerous for the resulting projective linear system to have dimension at least $4$.

Let
\begin{equation}\label{Fanod}
\mu=\mu_d:\p^5\dashrightarrow W\subseteq\p^N
\end{equation}
be a rational map defined by forms of degree $d\geq 2$ having multiplicity at least $d-1$ along $P$ and vanishing along $S$, and suppose that $\dim(W)=4$. In the terminology of the package \emph{SpecialFanoFourfolds} \cite{SpecialFanoFourfoldsSource}, $\mu$ is a \emph{Fano map} of $P\cup S$, and $W$ is the corresponding \emph{mirror fourfold}. The restriction to a secant line $\langle p,y\rangle$ of a form defining $\mu$ vanishes at $p$ with multiplicity at least $d-1$ and at $y$; hence it is determined up to a scalar, unless it vanishes identically. Thus, $\mu$ contracts every secant line to $P\cup S$ not contained in its base locus, and the general fibre of $\mu$ is a union of secant lines. If $X\subset\p^5$ is a cubic fourfold containing $P\cup S$, a general secant line meets $X$ in exactly one point outside $P\cup S$. Therefore, \( \mu|_X:X\dashrightarrow W \) is dominant, and its degree is the number of secant lines in a general fibre of $\mu$. In particular, $\mu|_X$ is birational if and only if the general fibre of $\mu$ is a single secant line; this condition does not depend on $X$, and it is satisfied in all our examples.

\begin{Proposition}\label{muQi}
Assume the hypotheses and notation of Proposition \ref{doublecoverT}, and let $\mu$ be a Fano map of $P\cup S$. Then, for every $i=1,\ldots,r$, the following hold.
\begin{enumerate}
\item[{\rm (i)}] If $e_i=0$ and $\langle Q_i\rangle\not\subset\Bs(\mu)$, then \( L_i=\overline{\mu(Q_i)}=\overline{\mu(\langle Q_i\rangle)}\subset W \) is a line or a point.

\item[{\rm (ii)}] If $e_i=1$ and $\langle Q_i\rangle\not\subset\Bs(\mu)$, then \( \overline{\mu(Q_i)}=\overline{\mu(\langle Q_i\rangle)} \) is a point $u_i\in W$.

\item[{\rm (iii)}] If $e_i\geq 2$, then $\langle Q_i\rangle\subset\Bs(\mu)$.
\end{enumerate}
\end{Proposition}

\begin{proof}
The restriction of $\mu$ to $\langle Q_i\rangle\simeq\p^3$ is defined by forms of degree $d$ having multiplicity at least $d-1$ along the plane $P\subset\langle Q_i\rangle$ and vanishing along $C_i\subset Q_i$.

If $e_i=0$, then $C_i$ is a line, and after removing the fixed component $(d-1)P$, the restricted linear system is contained in the pencil of planes containing $C_i$. Thus, $\mu$ maps $\langle Q_i\rangle\setminus\Bs(\mu)$ into a line. Since every plane containing $C_i$ meets $Q_i$ along $C_i$ and a second line, the closures of $\mu(Q_i)$ and $\mu(\langle Q_i\rangle)$ coincide.

If $e_i=1$, then $C_i$ is a conic spanning a plane. After removing the fixed component, the restricted linear system is spanned by an equation of the plane $\langle C_i\rangle$, and therefore $\mu$ contracts $\langle Q_i\rangle\setminus\Bs(\mu)$ to a point.

Finally, if $e_i\geq 2$, the curve $C_i\subset Q_i\subset\langle Q_i\rangle$ is nondegenerate. Hence every form in the linear system vanishes identically on $\langle Q_i\rangle$, so that \( \langle Q_i\rangle\subset\Bs(\mu). \)
\end{proof}

In particular, $L_i$ and $u_i$ do not depend on the cubic fourfold $X$ containing $P\cup S$, and $Q_i\subset\Bs(\mu|_X)$ if $e_i\geq 2$. In all our examples, $\langle Q_i\rangle\not\subset\Bs(\mu)$ whenever $e_i\leq 1$, and $L_i$ is a line whenever $e_i=0$.

We now define a key surface in $W$.

\begin{defin}\label{defU}
Let $(X,P)$ be a general pair, let $S\subset X$ be a section of $\pi_P|_X$, and let $\mu:\p^5\dashrightarrow W\subset\p^N$ be a Fano map of $P\cup S$. The surface
\[
U=\overline{\mu(E_S)}\subset W,
\]
where $E_S\subset X$ is the divisor swept out by the secant lines to $P\cup S$ contained in $X$, is called the \emph{surface associated with $(X,P,S,\mu)$}.
\end{defin}

\begin{rmk}\label{UeU'}{\rm
Since the secant lines to $P\cup S$ contained in $X$ are contracted by $\mu$, the closure $U$ of $\mu(E_S)$ is a surface. If the general fibre of $\mu$ is a single secant line, then $\mu$ induces a birational map \( \widetilde\mu:U_S\dashrightarrow U \) that sends a general secant line $[\ell]\in U_S$ to the point $\mu(\ell)\in U$. In this case, $U$ is birational to $K_X$; when $S$ is simple, we also denote by $\widetilde\mu$ the composition of $\widetilde\mu$ with the normalization morphism $\overline U_S\to U_S$. The surface $U$ depends on the degree $d$ of $\mu$.

If a map $\mu_d$ as above exists, then $\mu_{d+1}$ also maps $\p^5$ onto an irreducible fourfold whose general fibre is a secant line to $P\cup S$. Indeed, the hyperplanes containing $P$ define a base-point-free linear system on $\Bl_P(\p^5)$, and hence also on $\Bl_{P\cup S}(\p^5)$. As $d$ increases, the computations become more difficult, but the surface $U$ usually acquires better geometric properties.

Assume now that $S$ is simple. Fix $P\cup S$, and let $X'\subset\p^5$ be another general cubic fourfold containing it. If $U'\subset W$ is the surface associated with $(X',P,S,\mu)$, then \( \overline{\mu(\langle Q_i\rangle)}\subset U\cap U' \) whenever $e_i\leq 1$, since the quadric $Q_i$, as well as the corresponding quadric of $X'$, is contained in the divisor swept out by the secant lines. Thus, the lines $L_i$ and the points $u_i$ can be determined in two ways: either as the images under $\mu$ of the corresponding quadrics $Q_i$, or equivalently of their linear spans $\langle Q_i\rangle$, or from the intersection of the surfaces $U$ and $U'$ associated with two different general cubic fourfolds containing $P\cup S$.

Suppose now that $e_i\geq 2$. By Proposition \ref{muQi}, the birational map induced by $\mu|_X$ on $\Bl_{P\cup S}(X)$ has $Q_i$ in its base locus. Resolving the indeterminacies of this map, or equivalently performing suitable flops of the blown-up quadrics, produces further irreducible two-dimensional components in the base locus of $(\mu|_X)^{-1}$. Accordingly, the intersection $U\cap U'$ also contains curves corresponding to the indices $i$ with $e_i\geq 2$, which are described through the map $\widetilde\mu$ in Remark \ref{eigeq2} below. 
}
\end{rmk}

\begin{Proposition}\label{UStoU}
Assume the hypotheses and notation of Propositions \ref{doublecoverT} and \ref{muQi}, and suppose that the general fibre of $\mu$ is a single secant line to $P\cup S$. Then \( \widetilde\mu:U_S\dashrightarrow U \) is defined at every point $[\ell]\in U_S$ such that $\ell\not\subset\Bs(\mu)$, where it takes the value $\mu(\ell\setminus\Bs(\mu))$. Moreover, for every $i=1,\ldots,r$:
\begin{enumerate}
\item[{\rm (i)}] if $e_i=0$ and $L_i$ is a line, then $\widetilde\mu$ restricts to an isomorphism from $E_i$ onto $L_i$;

\item[{\rm (ii)}] if $e_i=1$ and $\langle Q_i\rangle\not\subset\Bs(\mu)$, then $\widetilde\mu$ contracts $E_i$ and $E'_i$ to the point $u_i$.
\end{enumerate}
\end{Proposition}

\begin{proof}
Let $[\ell_0]\in U_S$ be such that $\ell_0\not\subset\Bs(\mu)$. Then $[\ell_0]$ does not lie on a curve $\beta_X^{-1}(s'_i)$ with $e_i\geq 2$, because the lines parametrized by such a curve are contained in $\langle Q_i\rangle\subset\Bs(\mu)$ by Proposition \ref{muQi}. Hence $U_S$ is smooth at $[\ell_0]$, as observed after Proposition \ref{doublecoverT}. Let $x_0\in\ell_0\setminus\Bs(\mu)$. In a neighbourhood of $[\ell_0]$, in the analytic topology, the universal family of lines admits a section $\sigma$ with $\sigma([\ell_0])=x_0$. The map $\mu\circ\sigma$ is holomorphic near $[\ell_0]$, and it coincides with $\widetilde\mu$ on a dense open subset of $U_S$, because $\mu$ is constant on every secant line not contained in $\Bs(\mu)$. Hence $\widetilde\mu$ is defined at $[\ell_0]$, with value $\mu(x_0)$.

If $e_i=0$ and $L_i$ is a line, the proof of Proposition \ref{muQi} shows that $\Bs(\mu)\cap\langle Q_i\rangle=P\cup C_i$. A line $\ell$ parametrized by $E_i$ meets the line $C_i$ at one point $x$, and if $\ell\not\subset P$, then $\widetilde\mu$ maps $[\ell]$ to the point of $L_i$ corresponding to the plane $\langle\ell,C_i\rangle$, which is the tangent plane to $Q_i$ at $x$. Since $[\ell]\mapsto x$ is an isomorphism from $E_i$ onto $C_i$, and the tangent planes to $Q_i$ at the points of $C_i$ form the pencil of planes containing $C_i$, the restriction of $\widetilde\mu$ to $E_i$ is an isomorphism onto $L_i$. If $e_i=1$ and $\langle Q_i\rangle\not\subset\Bs(\mu)$, then $\Bs(\mu)\cap\langle Q_i\rangle=P\cup\langle C_i\rangle$, and the plane $\langle C_i\rangle$ contains no line of $Q_i$. Thus, $\widetilde\mu$ is defined at the general points of $E_i$ and $E'_i$, and it is constant there by Proposition \ref{muQi}.
\end{proof}

\begin{rmk}\label{eigeq2}{\rm
If $e_i\geq 2$, then $\langle Q_i\rangle\subset\Bs(\mu)$ by Proposition \ref{muQi}, and the behaviour of $\widetilde\mu$ along $E_i\cup E'_i$ is not determined by the restriction of $\mu$ to $\langle Q_i\rangle$. In our examples, the image of $E_i$ is a curve $\widetilde E_i\subset U$, which can be detected in $U\cap U'$ as explained in Remark \ref{UeU'}. In Examples 8, 10, 14 and 16 of Table \ref{TABex}, where $e_i=2$, the curve $\widetilde E_i$ is a conic, while in Example \ref{ex5} below, where $e_i=3$, it is a line consisting of non-normal points of $U$.
}
\end{rmk}

The following criterion identifies the projective models of $K_X$ obtained from $U$.

\begin{Proposition}\label{US=U}
Let $(X,P)$ be a general pair, let $S\subset X$ be a section of $\pi_P|_X$, and let $\mu$ be a Fano map of $P\cup S$ whose general fibre is a single secant line. Let $V$ be the blow-up of $K_X$ at finitely many distinct points, with exceptional curves $\Gamma_1,\ldots,\Gamma_m$, such that $\widetilde\mu$ induces a morphism $V\to U$, and let $H$ be the pullback to $V$ of a hyperplane section of $U$. Suppose that, for some integers $a_1,\ldots,a_m$, the complete linear system \( |H+a_1\Gamma_1+\cdots+a_m\Gamma_m| \) maps $V$ birationally onto a smooth surface $\widetilde U$ of sectional genus $\deg(\widetilde U)/2+1$. Then $\widetilde U$ is a K3 surface isomorphic to $K_X$, and $a_l=H\cdot\Gamma_l$ for every $l$. In particular,
\[
\deg(\widetilde U)
=
\deg(U)+\sum_{l=1}^m(H\cdot\Gamma_l)^2.
\]
\end{Proposition}

\begin{proof}
Let $\tau:V\to K_X$ be the blow-up morphism, let \( D=H+a_1\Gamma_1+\cdots+a_m\Gamma_m, \) and let $f:V\to\widetilde U$ be the morphism defined by $|D|$. Since $\widetilde U$ is smooth and birational to the K3 surface $K_X$, it is obtained from $K_X$ by a finite sequence of blow-ups at points, so that its canonical class $K_{\widetilde U}$ is effective. If $H_{\widetilde U}$ is a hyperplane section of $\widetilde U$, the adjunction formula gives \( K_{\widetilde U}\cdot H_{\widetilde U}=2g(\widetilde U)-2-\deg(\widetilde U)=0. \) Since $H_{\widetilde U}$ is ample, it follows that $K_{\widetilde U}=0$, and hence that $\widetilde U$ is a minimal surface. As birational maps between minimal surfaces of nonnegative Kodaira dimension are isomorphisms, $\sigma=f\circ\tau^{-1}$ is an isomorphism, and $f=\sigma\circ\tau$. Therefore \( D=f^*H_{\widetilde U}=\tau^*(\sigma^*H_{\widetilde U}) \) has degree zero on every curve $\Gamma_l$, that is, $a_l=H\cdot\Gamma_l$. Finally, \( \deg(\widetilde U)=D^2=H^2+\sum_l a_l^2, \) and $H^2=\deg(U)$ because $\widetilde\mu$ is birational.
\end{proof}

\begin{rmk}\label{dandg}
{\rm
In all our examples, the map $\widetilde\mu$ contracts the curves $E'_i$, and $H\cdot E_i$ is equal to $1$, $0$ or $e_i$ according to whether $e_i=0$, $e_i=1$ or $e_i\geq 2$; for $e_i\leq 1$, this agrees with Proposition \ref{UStoU}. Accordingly, Proposition \ref{US=U} is applied with $V=\overline U_S$ and the linear system
\[
\Bigl|H+\sum_{e_i=0}E_i+\sum_{e_i\geq 2}e_iE_i\Bigr|.
\]
Let $k$ be the number of indices $i$ with $e_i=0$, and let $g$ be the genus of a general curve in $|H|$, which is the sectional genus of $U$ when the singular locus of $U$ is finite. Then the degree and the sectional genus of $\widetilde U$ are
\[
\deg(\widetilde U)
=
\deg(U)+k+\sum_{e_i\geq 2}e_i^2
\qquad\text{and}\qquad
\widetilde g
=
g+\sum_{e_i\geq 2}\frac{e_i(e_i-1)}{2},
\]
and $\widetilde U\subset\p^{\widetilde g}$ is embedded by a complete linear system. If, in addition, the hyperplane sections of $U\subset\p^N$ cut out the complete linear system $|H|$ and $h^1(H)=0$, then the Riemann--Roch theorem gives \( N=g-k-\sum_{e_i\geq 2}e_i. \)

In some cases, the surface $U\subset\p^N$ may have finitely many singular points. In explicit computations, one must then either normalize $U$ or use, at least in principle, a Fano map $\mu_{d'}$ with $d'>d$. If $e_i=1$, the disjoint curves $E_i$ and $E'_i$ are contracted to the same point $u_i$, which is typically a non-normal point of $U$.

The Euler characteristic $\chi(\mathcal O_U)$ provides a useful consistency check. If $U$ has only rational singularities, then $\chi(\mathcal O_U)=\chi(\mathcal O_{K_X})=2$, whereas isolated non-normal points decrease $\chi(\mathcal O_U)$. For instance, $\chi(\mathcal O_U)=2-j$, where $j$ is the number of indices with $e_i=1$, in all the examples of Table \ref{TABex} except Examples 5 and 11.
}
\end{rmk}

We now apply these results to obtain geometric, and whenever possible explicit, projective realizations of the associated K3 surface $K_X$.

\begin{ex}\label{ex40}
We consider Example 18 in Table \ref{TABex}, corresponding to
\[
\mathfrak{S}\left(
\begin{array}{ccccc}
4 & 5 & 0 & 0 & 0\\
3 & 5 & 0 & 0 & 0
\end{array}
\right),
\]
which is the surface $S_{4,5}$ in the notation of Example \ref{Yer}; see also Section \ref{S53}. The surface $S\subset\p^5$ has degree $11$ and sectional genus $3$. The image of a general plane cubic through $p_1,\ldots,p_5$ is projected onto the curve $C=P\cap S$ of degree $7$ and geometric genus $1$, with $14$ nodes. Moreover, $\Sing(S)=\Sing(C)$, and $S$ is a simple section of the projection from $P$. We also have $p_i\notin\Delta_X$ for every $i=1,\ldots,5$.

The ideal of the reducible OADP surface $P\cup S\subset\p^5$ is generated by six cubic and three quartic forms. As shown in Section \ref{S53}, a general cubic fourfold $X$ containing $P\cup S$ is smooth, and
\[
[X]\in\mathcal S_{53}\subset\mathcal C_8\cap\mathcal C_{86}.
\]

The smallest degree of a Fano map of $P\cup S$ whose general fibre is a single secant line is $d=4$. The resulting map \( \mu_4:\p^5\dashrightarrow W\subset\p^{13} \) has image an irreducible fourfold $W$ of degree $20$, whose ideal is generated by $35$ quadratic forms. Its general fibre is a secant line to $P\cup S$, so that $\mu_4|_X$ is birational for every cubic fourfold $X$ containing $P\cup S$. The associated surface $U\subset W\subset\p^{13}$ has degree $29$, sectional genus $18$, and $\chi(\mathcal O_U)=2$.

The map $\pi_P|_S:S\dashrightarrow\p^2$ contracts only the images of the five exceptional curves on $\Bl_{p_1,\ldots,p_5}(\p^2)$, which are lines. Thus, $e_i=0$ for every $i=1,\ldots,5$, and $U_S$ is the blow-up of $K_X$ at five points by Proposition \ref{doublecoverT}. The surface $U\subset\p^{13}$ contains the five lines $L_i=\widetilde\mu(E_i)$, and one checks that it is smooth; in particular, $U$ is isomorphic to $U_S$.

The complete linear system \( |H+L_1+\cdots+L_5|, \) where $H$ is a hyperplane section of $U$, maps $U$ onto a smooth surface \( \widetilde U\subset\p^{18} \) of degree $34$ and sectional genus $18$, whose ideal is generated by $120$ quadratic forms. By Proposition \ref{US=U}, $\widetilde U$ is a K3 surface isomorphic to $K_X$.

Let $H_X$ be the hyperplane class of the embedding $K_X\simeq\widetilde U\subset\p^{18}$, and let $D_X$ be the degree-$2$ polarization defining the double cover \( \alpha_X:K_X\longrightarrow\p^2. \) By computing the inverse of $\mu_4|_X$, we find that the lattice $\langle H_X,D_X\rangle$ has discriminant $-\delta=-53$; see also \cite{HassettPlane}. Therefore,
\[
(H_X\cdot D_X)^2
=
53+D_X^2H_X^2
=
53+2\cdot 34
=
121,
\]
and hence \( H_X\cdot D_X=11. \) This equality can also be verified by pulling back to $U$ the strict transform on $S$ of a general line in $\p^2$, using the inverses of $\mu_4|_X$ and $\pi_P|_S$, together with Proposition \ref{doublecoverT}.
\end{ex}

\begin{ex}\label{ex17}
We consider Example 7 in Table \ref{TABex}, corresponding to
\[
\mathfrak{S}\left(
\begin{array}{ccccc}
4 & 4 & 1 & 0 & 0\\
1 & 0 & 0 & 0 & 0
\end{array}
\right).
\]
The surface $S\subset\p^5$ has degree $8$ and sectional genus $2$. The line $E$ is projected onto the plane quartic $C=P\cap S$, whose three nodes constitute the entire singular locus of $S$. The surface $S$ is a simple section of the projection from $P$, and the base points of the inverse of $\pi_P|_S$ do not lie on $\Delta_X$.

The ideal of the reducible OADP surface $P\cup S\subset\p^5$ is generated by one quadratic and seven cubic forms. A general cubic fourfold $X$ containing $P\cup S$ is smooth, and formula \eqref{deltaex} gives $\delta=45$. Formula \eqref{selfS} also shows that
\[
[X]\in\mathcal C_8\cap\mathcal C_{32}.
\]

Since $P\cup S$ is contained in a unique quadric, there is no Fano map of degree $2$, and the smallest degree of a Fano map is $d=3$: the cubic forms that are singular along $P$ and vanish along $S$ define a map \( \mu_3:\p^5\dashrightarrow W\subset\p^6, \) whose image is an irreducible complete intersection of two quadrics. Its general fibre is a secant line to $P\cup S$, and its restriction to a general cubic fourfold $X$ containing $P\cup S$ is birational. The associated surface $U\subset W\subset\p^6$ has degree $14$, sectional genus $10$, and $\chi(\mathcal O_U)=1$.

The map $\pi_P|_S:S\dashrightarrow\p^2$ contracts five curves. Four of them are the images of the strict transforms of the lines $\langle p_i,p_5\rangle$, for $i=1,\ldots,4$, while the fifth is the image of the strict transform of the unique conic through $p_1,\ldots,p_5$. Therefore,
\[
e_i=0\quad\text{for }i=1,\ldots,4,
\qquad
e_5=1.
\]
By Proposition \ref{doublecoverT}, $U_S$ is the blow-up of $K_X$ at the six points $s_1,\ldots,s_5,s'_5$. The surface $U\subset\p^6$ is non-normal, with a unique singular point \( u_5=\mu_3(Q_5), \) to which $\widetilde\mu$ contracts the curves $E_5$ and $E'_5$, and it contains four lines $L_1,\ldots,L_4$, none of which passes through $u_5$.

Let $U_1\subset\p^5$ be the projection of $U$ from $u_5$. Then $U_1$ is smooth and isomorphic to $U_S$. Let $H_1$ be the pullback to $U_1$ of the hyperplane class of $U$, and let $L'_1,\ldots,L'_4\subset U_1$ be the strict transforms of $L_1,\ldots,L_4$. The complete linear system \( |H_1+L'_1+\cdots+L'_4| \) maps $U_1$ onto a surface $\widetilde U\subset\p^{10}$ of degree $18$ and sectional genus $10$, whose homogeneous ideal is generated by $28$ quadratic forms. The surface $\widetilde U$ is smooth, and hence isomorphic to $K_X$ by Proposition \ref{US=U}.

Let $H_X$ and $D_X$ be the polarizations on $K_X$ defined as in the previous example. Since the lattice $\langle H_X,D_X\rangle$ has discriminant $-45$, we have
\[
(H_X\cdot D_X)^2
=
45+D_X^2H_X^2
=
45+2\cdot18
=
81.
\]
Therefore, \( H_X\cdot D_X=9. \) As in Example \ref{ex40}, this intersection number can also be verified directly by pulling back a general line in $\p^2$ through the relevant birational maps.
\end{ex}

In some cases, one may have $\deg(\widetilde E_i)<e_i$ for some $e_i\geq 2$. The curve $\widetilde E_i\subset U$ then consists of singular points of $U$, as illustrated by the following example.

\begin{ex}\label{ex5}
We briefly consider the following example:
\[
\mathfrak{S}\left(
\begin{array}{ccccc}
6 & 10 & 0 & 0 & 1\\
1 & 2 & 0 & 0 & 0
\end{array}
\right).
\]
The surface $S\subset\p^5$ has degree $10$ and sectional genus $4$. The line $E$ through two of the ten simple base points is projected onto the plane quartic $C=P\cap S$, whose three nodes constitute the entire singular locus of $S$. The surface $S$ is a simple section of the projection from $P$.

The ideal of the reducible OADP surface $P\cup S\subset\p^5$ is generated by seven cubic and three quartic forms. A general cubic fourfold $X$ containing $P\cup S$ is smooth, and formula \eqref{deltaex} gives $\delta=21$. Moreover,
\[
[X]\in\mathcal C_8\cap\mathcal C_{14}.
\]

The smallest degree of a Fano map of $P\cup S$ whose general fibre is a single secant line is $d=4$. The resulting map \( \mu_4:\p^5\dashrightarrow W\subset\p^{10} \) has image a fourfold $W$ of degree $10$, whose ideal is generated by $18$ quadratic forms. Its restriction to a general cubic fourfold $X$ containing $P\cup S$ produces a surface $U\subset W\subset\p^{10}$ of degree $31$, sectional genus $25$, and $\chi(\mathcal O_U)=2$.

The map $\pi_P|_S:S\dashrightarrow\p^2$ contracts ten lines and a smooth rational quartic curve. Hence
\[
e_i=0\quad\text{for }i=1,\ldots,10,
\qquad
e_{11}=3.
\]
By Proposition \ref{doublecoverT}, $U_S$ is not normal along $\beta_X^{-1}(s'_{11})$, and its normalization $\overline U_S$ is the blow-up of $K_X$ at the twelve points $s_1,\ldots,s_{11},s'_{11}$. The surface $U\subset\p^{10}$ contains the ten lines $L_1,\ldots,L_{10}$. A direct computation with \emph{Macaulay2} shows that the intersection of the surfaces associated with two general cubic fourfolds containing $P\cup S$ consists of $11$ lines; hence $\widetilde E_{11}=\widetilde\mu(E_{11})$ is a line. We shall see below that $E_{11}$ is mapped onto $\widetilde E_{11}$ with degree $3$, so that $\widetilde E_{11}$ consists of non-normal points of $U$ of multiplicity $3$. In particular, $U$ is not normal along a curve, although $\chi(\mathcal O_U)=2$; compare Remark \ref{dandg}.

Let $H'$ be the pullback to $\overline U_S$ of the hyperplane class of $U$. By Bertini's theorem, a general curve in $|H'|$ is smooth, has degree $31$, and has genus $23$. The complete linear system \( |H'+E_1+\cdots+E_{10}+3E_{11}| \) maps $\overline U_S$ birationally onto a smooth surface $\widetilde U\subset\p^{26}$ of degree \( 31+10+9=50 \) and sectional genus \( 23+3=26, \) whose ideal is generated by $276$ quadratic forms. By Proposition \ref{US=U}, applied to the blow-up $\overline U_S\to K_X$, the surface $\widetilde U$ is a K3 surface isomorphic to $K_X$; moreover, $H'\cdot E_{11}=3$ and $H'\cdot E'_{11}=0$. Thus, $E_{11}$ is mapped onto the line $\widetilde E_{11}$ with degree $3$, while $E'_{11}$ is contracted by $\widetilde\mu$.

The construction and singularities of $U$ become clearer when this procedure is reversed. Imposing a triple point at the image of $E_{11}$ on the linear system of hyperplane sections produces a surface \( \widehat U\subset\p^{20} \) containing a twisted cubic $\widehat C_{11}$. To recover $U$, one must choose ten points on $\widehat U$ whose linear span is a $\p^9$ meeting $\langle\widehat C_{11}\rangle$ along a line. By contrast, the span of ten general points of $\widehat U$ is disjoint from $\langle\widehat C_{11}\rangle$. Thus, $U$ is obtained from $\widetilde U$ by a special projection.

An explicit construction of $\widetilde U$ can be obtained as follows. Projecting $U\subset\p^{10}$ from the line $\widetilde E_{11}$ gives a smooth surface \( U_1\subset\p^8 \) of degree $24$ and sectional genus $20$. The projection induces a birational morphism \( U_1\longrightarrow U, \) and the strict transform of $\widetilde E_{11}$ on $U_1$ is a smooth rational normal quartic curve \( C_{11}\subset U_1. \) It follows that $U_1\simeq\overline U_S$ and that $E_{11}$ is mapped isomorphically onto $C_{11}$.

Let $H_1$ be the hyperplane class of $U_1$, and let $L'_1,\ldots,L'_{10}\subset U_1$ be the strict transforms of the lines $L_1,\ldots,L_{10}$. The complete linear system \( |H_1+L'_1+\cdots+L'_{10}+4C_{11}| \) maps $U_1$ birationally onto a smooth K3 surface \( \widetilde U_1\subset\p^{26} \) of degree \( 24+10+16=50 \) and sectional genus \( 20+6=26. \) The surface $\widetilde U_1$ is projectively equivalent to $\widetilde U$. In this case, the conditions imposed by the quadruple curve and the ten simple curves are not independent, since \( 10+10+8=28>26. \) Let $H_X$ be the polarization defining the embedding $K_X\simeq\widetilde U\subset\p^{26}$, and let $D_X$ be the degree-$2$ polarization defining the double cover of $\p^2$. Since the lattice $\langle H_X,D_X\rangle$ has discriminant $-21$, we have
\[
(H_X\cdot D_X)^2
=
21+D_X^2H_X^2
=
21+2\cdot50
=
121.
\]
Therefore, \( H_X\cdot D_X=11. \)
\end{ex}

The examples listed in Table \ref{TABex} can be treated by the same methods, producing explicit projective models $K_X\simeq\widetilde U\subset\p^{\widetilde g}$ together with the lattices $\langle H_X,D_X\rangle$; see Section \ref{Computations} for a computational illustration.

\section{Computations with \emph{Macaulay2}}\label{Computations}

The examples and computations presented in this paper can be studied using the
\emph{Macaulay2} package \emph{SpecialFanoFourfolds}
\cite{SpecialFanoFourfoldsSource,2024packageSFF,macaulay2}. We illustrate some of its
functionalities using Example \ref{ex40}, which corresponds to a surface of
type
\[
\mathfrak{S}\left(
\begin{array}{cccc}
4 & 5 & 0 & \cdots\\
3 & 5 & 0 & \cdots
\end{array}
\right).
\]
The computations below should be run with version 1.26.11 of
\emph{Macaulay2}.
The mathematical output displayed here is obtained in
\texttt{WebApp} mode, that is, with
\texttt{topLevelMode = WebApp}; in this mode, the package provides a native
\LaTeX{} representation. Some outputs have been abbreviated and slightly
edited for print.

The following commands load the package and construct a random cubic
fourfold containing the corresponding surface and plane.

\begin{tcolorbox}[
  colback=gray!4,
  colframe=gray!55,
  boxrule=0.4pt,
  arc=1mm,
  left=2mm,
  right=2mm,
  top=2mm,
  bottom=2mm
]
{\footnotesize\ttfamily
i1 : needsPackage "SpecialFanoFourfolds"; \hfill -- version 2.8.1

\medskip

i2 : X = cubicFourfold surface((4,5,0,0,0),(3,5,0,0,0));
}

\medskip

\noindent
{\footnotesize\ttfamily o2 :}\quad
{\footnotesize\ttfamily ProjectiveVariety, }\(
\text{cubic fourfold in }\mathcal C_{86}\cap\mathcal C_8
\)
\end{tcolorbox}

The command \texttt{describe} summarizes the main geometric and numerical
properties of the example, while \texttt{surfaces} returns the surface $S$
and the plane $P$.

\begin{tcolorbox}[
  colback=gray!4,
  colframe=gray!55,
  boxrule=0.4pt,
  arc=1mm,
  left=2mm,
  right=2mm,
  top=2mm,
  bottom=2mm
]
\noindent
{\footnotesize\ttfamily i3 : describe X}

\medskip

\noindent
{\footnotesize\ttfamily o3 =}\quad
\begin{minipage}[t]{0.9\linewidth}
\(
\begin{aligned}[t]
&\text{Cubic fourfold in }\mathcal C_{86}\cap\mathcal C_8
  \text{ of lattice discriminant}
\\
&\det\left(
\begin{array}{ccc}
3  & 11 & 1\\
11 & 69 & -4\\
1  & -4 & 3
\end{array}
\right)
=53=6\cdot 8+5,
\\
&\text{containing two surfaces:}
\\
&\quad
  \text{- a 14-nodal rational surface of degree }11
  \text{ and sectional genus }3,
\\
&\quad
  \text{- a plane},
\\
&\text{whose intersection is a 14-nodal curve of degree }7
  \text{ and arithmetic genus }15,
\\
&\text{K3 status: [░░░░░ / ▓▓▓▓▓]\qquad \{ID: 40\}.}
\end{aligned}
\)
\end{minipage}

\medskip

\noindent
{\footnotesize\ttfamily i4 : (S,P) = surfaces X}

\medskip

\noindent
{\footnotesize\ttfamily o4 =}\quad
\(\left(S,P\right)\)
\end{tcolorbox}

Thus, $X$ belongs to $\mathcal C_{86}\cap\mathcal C_8$, and the lattice
$\langle h^2,S,P\rangle$ has discriminant $53$. The surface $S$ has degree
$11$ and sectional genus $3$, while $C=S\cap P$ is a plane curve of degree
$7$ with arithmetic genus $15$ and geometric genus $1$. Its $14$ nodes
constitute the singular locus of both $C$ and $S$.

The package also implements the parameter count discussed in Proposition
\ref{moduliSC} and Section \ref{ParCount}. Setting
\texttt{Verbose=>true} displays the main steps of the computation; some
output lines are omitted below.

{\footnotesize
\begin{M2Verbatim}
i5 : parameterCount(X,Verbose=>true)

-- starting parameterCount computation
-- input: cubic fourfold X containing two surfaces:
   -- S = surface((4,5,0,0,0),(3,5,0,0,0)):
      14-nodal rational surface of degree 11 and sectional genus 3
   -- (projected from PP^9 with center a PP^3 ⊂ PP^6 (⊂ PP^9))
   -- P: plane
   -- C = S ∩ P: curve in PP^2 defined by a form of degree 7

-- moduli count for 5 points in ℙ²: 2
-- dimension of the space of plane curves of degree 3
   passing through 5 general points: 4
-- dim GG(3,6) = 12
-- taking dim Aut(S,ℙ⁵) = 0
-- dimension of the family of reducible surfaces S ∪ P in ℙ⁵:
   2 + 4 + 12 + (35 - 0) = 53

-- h^0(I_{S ∪ P,ℙ⁵}(3)) = 6
-- dimension of the incidence variety {(S,P,X) : S ∪ P ⊂ X}:
   53 + 5 = 58

-- computing the normal sheaf of S in X...
-- h^0(N_{S,X}) = 5
-- computing the normal sheaf of P in X...
-- h^0(N_{P,X}) = 0
-- upper bound for the dimension of the family of reducible
   surfaces S ∪ P in X: 5

-- codim. in C_8 of {[X] : S ∪ P ⊂ X}
   ≤ 54 - (58 - 5) = 1 ✅

o5 = (1,(6,53,5))
\end{M2Verbatim}
}

The output records the codimension estimate together with the three integers
\[
\bigl(
\h^0(\mathcal I_{S\cup P}(3)),
n(\mathfrak S),
\h^0(N_{S/X})+\h^0(N_{P/X})
\bigr)
=
(6,53,5).
\]
Here $\h^0(N_{P/X})=0$, as expected since $X$ is smooth.
Therefore
\[
53+6-5=54,
\]
exactly as required in Proposition \ref{suffSd}. The computation gives an
upper bound of $1$ for the codimension of $\mathcal C_{\mathfrak S}$ in
$\mathcal C_8$, and Proposition \ref{suffSd} identifies $\mathcal C_{\mathfrak S}$
with $\mathcal S_{53}$; this is the computation used in Section \ref{S53}.
A refinement of this parameter count is described in Remark
\ref{refinedCount}; for the example above, it gives the same bound. It will
be used by \texttt{parameterCount} in a future version of the package, which
can be installed with the function \texttt{updatePackage} before it is
distributed with \emph{Macaulay2}.

We now construct the K3 surface associated with $X$. On its first call, the
command \texttt{polarizedK3surface} computes the underlying projective K3
surface. This construction involves computing the Fano map $\mu$, extracting a
non-minimal surface $U$, identifying its exceptional curves, and contracting
them to obtain the minimal K3 surface. The option \texttt{Strategy} controls
the methods used in this construction and in the subsequent computation of
the lattice polarization. Here the first stage uses the strategy
\texttt{"Inverse"}, which computes the inverse of $\mu|_X$, and the
polarization uses the strategy \texttt{"Genus2Curve"}, which pulls back to
$U$ the strict transform on $S$ of a general line in $\mathbb P^2$, as at the
end of Example \ref{ex40}. These are not the default choices, which are
faster in this example. Only the most relevant output lines are reproduced
below.

{\footnotesize
\begin{M2Verbatim}
i6 : E = polarizedK3surface(X,
         Strategy=>("Inverse","Genus2Curve"), Verbose=>true)

-- starting underlying K3 computation
-- planned steps:
-- 1. compute Fano map μ : ℙ⁵ ⇢ W
-- 2. extract surface U from the base locus of (μ|X)⁻¹
-- 3. extract U' from another general cubic X'
-- 4. analyze U ∩ U' for exceptional curves
-- 5. compute the contraction map f : U -> Ũ
-- 6. prepare data for lattice polarization

-- fanoMap: attempting maps of degrees 2, 3, and 4...
-- computed map μ: ℙ⁵ --> ℙ¹³ defined by hypersurfaces
   of degree 4 with multiplicities 1 and 3 along the surfaces
-- verified: generic fiber is a line connecting a point
   on each of the two surfaces
-- image of μ: 4-dimensional subvariety of ℙ¹³
   cut out by 35 hypersurfaces of degree 2
-- surface U in ℙ¹³ cut out by 50 hypersurfaces of degree 2
-- extracting the 1-dimensional part of U ∩ U'...
-- number of lines found: 5
-- computing the contraction map from U to the minimal K3 surface
\end{M2Verbatim}
}

\begin{tcolorbox}[
  colback=gray!4,
  colframe=gray!55,
  boxrule=0.4pt,
  arc=1mm,
  left=2mm,
  right=2mm,
  top=2mm,
  bottom=2mm
]
\noindent
{\footnotesize\ttfamily o6 =}\quad
\begin{minipage}[t]{0.9\linewidth}
\(
\begin{aligned}[t]
&\text{K3 surface associated to }X
\\
&\text{Mirror fourfold: a fourfold in }\mathbb P^{13}
  \text{ of degree }20
\\
&\qquad
  \text{and sectional genus }10,
  \text{ cut out by }35\text{ quadrics}
\\
&\text{Surface }U\text{: a surface of degree }29
  \text{ and sectional genus }18,
\\
&\qquad
  \chi(\mathcal O_U)=2,
  \text{ cut out by }50\text{ quadrics}
\\
&\text{Exceptional curves: }5\text{ lines}
\\
&\text{Minimal K3 surface }\widetilde U\text{: a surface in }
  \mathbb P^{18}
\\
&\qquad
  \text{of degree }34\text{ and sectional genus }18,
  \text{ cut out by }120\text{ quadrics}
\\
&\text{Lattice polarization: not yet computed;}
\\
&\qquad
\text{rerun }\texttt{polarizedK3surface}.
\end{aligned}
\)
\end{minipage}

\medskip

\noindent
{\footnotesize\ttfamily o6 :}\quad
Lattice-polarized K3 surface associated to \(X\)

{\scriptsize\ttfamily \qquad\qquad -- K3 status: [▓▓▓▓░ / ▓▓▓▓▓]}
\end{tcolorbox}

At this stage, the underlying K3 surface has been constructed. The output
summarizes the main geometric objects appearing in the construction: the
\emph{mirror fourfold} $W$, the non-minimal surface $U$, its five exceptional
lines, and the minimal K3 surface
$\widetilde U\simeq K_X\subset\mathbb P^{18}$
of degree $34$ and sectional genus $18$ (see Example \ref{ex40}). A second call to
\texttt{polarizedK3surface} computes the lattice polarization.

{\footnotesize
\begin{M2Verbatim}
i7 : E = polarizedK3surface(X, Verbose=>true)

-- starting polarization computation
-- settings: Verbose => true, Strategy => "Genus2Curve"
-- obtained a curve on U of degree 11 and genus 2
-- computing its image on the K3 surface...
-- constructing lattice-polarized K3 with (g,d,C^2) = (18,11,2)
-- polarization successfully completed
\end{M2Verbatim}
}

\begin{tcolorbox}[
  colback=gray!4,
  colframe=gray!55,
  boxrule=0.4pt,
  arc=1mm,
  left=2mm,
  right=2mm,
  top=2mm,
  bottom=2mm
]
\noindent
{\footnotesize\ttfamily o7 =}\quad
\begin{minipage}[t]{0.9\linewidth}
\(
\begin{aligned}[t]
&\text{K3 surface associated to }X
\\
&\qquad\vdots
\\
&\text{Lattice intersection matrix on }\widetilde U\text{: }
\left(
\begin{array}{cc}
34 & 11\\
11 & 2
\end{array}
\right)
\quad
\bigl(\det=-53\bigr)\ \checkmark
\\
&\bigl[
  \text{use }\texttt{E(a,b)}
  \text{ or }\texttt{E(g)}
  \text{ for polarized K3 models}
  \bigr]
\\
&\bigl[
  \text{call }\texttt{latticePolarization(E)}
  \text{ for the lattice polarization}
  \bigr].
\end{aligned}
\)
\end{minipage}

\medskip

\noindent
{\footnotesize\ttfamily o7 :}\quad
Lattice-polarized K3 surface associated to \(X\)

{\scriptsize\ttfamily \qquad\qquad -- K3 status: [▓▓▓▓▓ / ▓▓▓▓▓]}
\end{tcolorbox}

The object \texttt{E} contains the construction data and the lattice polarization.
The command \texttt{building E} provides access to intermediate data from
the construction, such as the Fano map and the surface $U$, while
\texttt{latticePolarization E} returns the corresponding rank-$2$ lattice
$\langle H,C\rangle$, where $H$ is the hyperplane class of
$\widetilde U\subset\mathbb P^{18}$ and $C$ is the class of the image of the
genus-$2$ curve obtained above; in the notation of Example \ref{ex40},
$H=H_X$ and $C=D_X$. Its output also suggests several divisor classes that
can be used to construct projective models of the K3 surface.

\begin{tcolorbox}[
  colback=gray!4,
  colframe=gray!55,
  boxrule=0.4pt,
  arc=1mm,
  left=2mm,
  right=2mm,
  top=2mm,
  bottom=2mm
]
\noindent 
{\footnotesize\ttfamily i8 : latticePolarization E} 

\medskip 
\noindent
{\footnotesize\ttfamily o8 =}\quad
\begin{minipage}[t]{0.9\linewidth}
\(
\begin{aligned}[t]
&\text{Lattice rank-2 polarization with intersection matrix }
\left(
\begin{array}{cc}
34&11\\
11&2
\end{array}
\right)
\\
&(a,b)=(0,1)\longrightarrow g=2
\\
&(a,b)=(1,-1)\longrightarrow g=8
\\
&(a,b)=(1,0)\longrightarrow g=18
\\
&(a,b)=(1,1)\longrightarrow g=30
\\
&(a,b)=(1,2)\longrightarrow g=44
\\
&\qquad\vdots
\\
&[\text{more lines with }\texttt{polarize(i,oo)}]
\\
&[\text{use }\texttt{map(oo,a,b)}
  \text{ to obtain the corresponding map}]
\\
&[\text{on }\texttt{surface(oo)}
  \text{ defined by }|aH+bC|,
  \text{ where }(H,C)=\texttt{basis(oo)}].
\end{aligned}
\)
\end{minipage}

\medskip

\noindent
{\footnotesize\ttfamily o8 :}\quad
Lattice-polarization on K3 surface associated to \(X\)
\end{tcolorbox}

More generally, \texttt{E(a,b)} constructs the projective model defined by
the divisor class $aH+bC$ in the current lattice basis. Alternatively,
\texttt{E(g)} searches for a suitable divisor class of genus $g$ and
constructs the corresponding model whenever one is found. The same commands may be applied
to any projective model obtained in this way.

For instance, the divisor class $H-C$, corresponding to $(a,b)=(1,-1)$, has
self-intersection $14$ and therefore genus $8$. The corresponding
projective model can be constructed directly from \texttt{E}.

{\footnotesize
\begin{M2Verbatim}
i9 : E8 = E(1,-1,Verbose=>true)

-- (*) constructing K3 surface of genus 8 in PP^8
   -- from K3 surface of genus 18 in PP^18
   -- using divisor D = aH+bC, with D^2 = 14 = 2*8-2,
      where (a,b) = (1,-1)
   -- computing K3 surface as image of map to PP^8
-- (*) constructing polarization on K3 surface of genus 8 in PP^8...
-- (✓) polarized K3 surface of genus 8 in PP^8 successfully constructed
\end{M2Verbatim}
}

\begin{tcolorbox}[
  colback=gray!4,
  colframe=gray!55,
  boxrule=0.4pt,
  arc=1mm,
  left=2mm,
  right=2mm,
  top=2mm,
  bottom=2mm
]
\noindent
{\footnotesize\ttfamily o9 =}\quad
\begin{minipage}[t]{0.9\linewidth}
\(
\begin{aligned}[t]
&\text{Lattice-polarized K3 surface of genus }8
  \text{ in }\mathbb P^8
\\
&\text{of lattice discriminant }
\det\left(
\begin{array}{cc}
14&9\\
9&2
\end{array}
\right)
=-53
\\
&\text{associated to a cubic fourfold in }
\mathcal C_{86}\cap\mathcal C_8
\text{ of lattice discriminant}
\\
&\det\left(
\begin{array}{ccc}
3&11&1\\
11&69&-4\\
1&-4&3
\end{array}
\right)
=53.
\end{aligned}
\)
\end{minipage}

\medskip

\noindent
{\footnotesize\ttfamily o9 :}\quad
Lattice-polarized K3 surface associated to \(X\)
\end{tcolorbox}

The resulting object \texttt{E8} represents a genus-$8$ projective model of the
same K3 surface. As with \texttt{E}, one may inspect its lattice polarization or use
it to construct further projective models. Its lattice basis consists of the
hyperplane class of the genus-$8$ model and the image of the curve $C$
above. In particular, the divisor class with coordinates $(0,1)$ has
self-intersection $2$ and yields the genus-$2$ model, namely the double cover
$\alpha_X\colon K_X\to\mathbb P^2$ branched along the discriminant sextic
$\Delta_X$, realized as a sextic surface in $\mathbb P(1,1,1,3)$.

{\footnotesize
\begin{M2Verbatim}
i10 : E2 = E8(0,1,Verbose=>true)

-- (*) constructing K3 surface of genus 2 in PP(1,1,1,3)
   -- from K3 surface of genus 8 in PP^8
   -- using divisor D = aH+bC, with D^2 = 2 = 2*2-2,
      where (a,b) = (0,1)
   -- obtained map from K3 surface in PP^8 to PP(1,1,1,3)
   -- computing the image in PP(1,1,1,3)...
-- (*) constructing polarization on K3 surface of genus 2
   in PP(1,1,1,3)...
-- (✓) polarized K3 surface of genus 2 in PP(1,1,1,3)
   successfully constructed
\end{M2Verbatim}
}

\begin{tcolorbox}[
  colback=gray!4,
  colframe=gray!55,
  boxrule=0.4pt,
  arc=1mm,
  left=2mm,
  right=2mm,
  top=2mm,
  bottom=2mm
]
\noindent
{\footnotesize\ttfamily o10 =}\quad
\begin{minipage}[t]{0.9\linewidth}
\(
\begin{aligned}[t]
&\text{Lattice-polarized K3 surface of genus }2
  \text{ in }\mathbb P(1,1,1,3)
\\
&\text{of lattice discriminant }
\det\left(
\begin{array}{cc}
2&9\\
9&14
\end{array}
\right)
=-53
\\
&\text{associated to a cubic fourfold in }
\mathcal C_{86}\cap\mathcal C_8
\text{ of lattice discriminant}
\\
&\det\left(
\begin{array}{ccc}
3&11&1\\
11&69&-4\\
1&-4&3
\end{array}
\right)
=53.
\end{aligned}
\)
\end{minipage}

\medskip

\noindent
{\footnotesize\ttfamily o10 :}\quad
Lattice-polarized K3 surface associated to \(X\)
\end{tcolorbox}

\section{Summary tables of examples}\label{SummaryTables}

Table \ref{TABex} presents a selection of explicit examples illustrating the
constructions considered in this paper. For each example, it records the
main construction data, the relevant Noether--Lefschetz loci, the parameter
count, and information on the associated polarized K3 surface. The table
was generated with the function \texttt{texTables} of the
\emph{SpecialFanoFourfolds} package.

The column \emph{parameter count} contains the upper bound for
$\codim_{\mathcal C_8}\mathcal C_{\mathfrak S}$ obtained in Section \ref{ParCount},
in the form $54-(n(\mathfrak S)+\h^0(\mathcal I_{S\cup P}(3))-1-\h^0(N_{S/X}))$.
The columns (\textdagger) and (\textdaggerdbl) record, respectively, the pairs
\[
\bigl(n(\mathfrak S),\h^0(N_{S/\mathbb P^5})\bigr)
\qquad\text{and}\qquad
\bigl(
\h^0(\mathcal I_{S\cup P}(3)),
\h^0(\mathcal I_{S}(3))
\bigr).
\]
The columns (\textsection) and (\textasteriskcentered) record, respectively,
\(\h^0(N_{S/X})\) and the upper bound for
\[
\codim_{\mathcal C}
\overline{\{[X]\in\mathcal C:S\subset X\}}
\]
given by the parameter count of \cite{Nuer-Unir}, as implemented in the
package. In the latter column, (!) indicates that one of the conditions
$\h^1(\mathcal O_S(3))=0$ and $\h^0(\mathcal I_S(3))=56-\chi(\mathcal O_S(3))$,
which the package checks before performing this count, is not satisfied, so
that the corresponding command stops with an error; the table reports the
value that the count would otherwise give. A dash indicates that the value has
not been computed. For a fourfold \texttt{X} constructed as in Section
\ref{Computations}, the count in the column (\textasteriskcentered) is obtained
with \texttt{parameterCount(new CubicFourfold from X)}, whereas
\texttt{parameterCount X} gives the count in $\mathcal C_8$.

In the column of Noether--Lefschetz loci, the first two divisors are
$\mathcal C_d$, with $d=\disc\langle h^2,S\rangle$ computed by means of
\eqref{selfS}, and $\mathcal C_8$. This list is not exhaustive, since $\mathcal S_\delta$ is
contained in infinitely many Hassett divisors; see Section \ref{SdeltaSec}.

In Examples 5, 8, 11, 14, 17, 19, 20, 21, 23 and 24 of Table \ref{TABex},
one has $\chi(\mathcal O_U)\neq 2$, and the construction of $\widetilde U$
requires the normalization of $U$, which is very expensive to compute. In
these cases, the invariants of $\widetilde U$ and its lattice polarization
are obtained \emph{virtually}, that is, from $U$ and from curves on $U$,
without constructing $\widetilde U$. This is only a matter of computation
time: with a non-virtual strategy of \texttt{polarizedK3surface}, the
computation can in principle be completed.

Some of the computations associated with these examples can be
time-consuming. The examples in Table \ref{TABex}, together with several
others, are available in an online archive maintained with the package.
After importing the archive,
a stored example of the same type as Example
\ref{ex40} can be retrieved by its identifier as follows.

{\footnotesize
\begin{M2Verbatim}
i11 : A = localExamplesArchive()

i12 : importFrom("@",A)

i13 : X = example(40,0);
\end{M2Verbatim}
}

For the current versions of the article and the archive, 
the identifiers corresponding to the rows of Table \ref{TABex} are 
\[ \begin{array}{c|rrrrrrrrrrrr} i&1&2&3&4&5&6&7&8&9&10&11&12\\ \mathrm{ID}&3&12&6&9&11&26&17&19&21&22&23&36\\[2pt] i&13&14&15&16&17&18&19&20&21&22&23&24\\ \mathrm{ID}&30&31&33&34&35&40&25&28&38&39&32&37 \end{array}. \] 
Thus, the example in row \(i\) can be retrieved with \texttt{example(ID,0)}.

The stored examples retain available data on parameter counts, associated K3
surfaces, and other computations provided by the package, avoiding the need
to repeat the most time-consuming calculations. The online documentation of
the package contains further details on the archive and on the functions
\texttt{localExamplesArchive} and \texttt{texTables}.

\newgeometry{
  left=0.5cm,
  right=0.5cm,
  top=0.5cm,
  bottom=0.5cm
}
\pagestyle{empty}
\thispagestyle{empty}
\refstepcounter{table}
\label{TABex}
\vspace*{-\topskip}
\noindent
\rotatebox{90}{%
\begin{adjustbox}{
  width=0.98\textheight,
  totalheight=0.98\textwidth,
  keepaspectratio
}

\begin{tabular}{|c||c|c|l|c||l|c|c|c|c||c|c|c|c|c|r||c|}
\hline
$i$ & surface $S\subset X\subset \mathbb{P}^5$ & $S, C=S\cap P$ & & & parameter count &(\textdagger) & (\textdaggerdbl) & (\textsection) & (\textasteriskcentered) & $\mu:\mathbb{P}^5\dashrightarrow W$ & $W$ & $U$ & exc. curves & $\widetilde{U}$ & & $i$ \\
\hline
\\
\hline
1 & $\mathfrak{S}\left(\!\begin{array}{cccc}
4&9&0&0\\
3&9&0&0
\end{array}\!\right) $ & $\begin{array}{c} \deg(S)=7 , g(S)= 3 \\ \operatorname{gens\,id.}: {2}^{3}\ {3}^{1}\  \\ \deg(C)= 3 , n(C)= 0 \end{array}$  & $21 = 2 \cdot 8 + 5 =  \det \left(\!\begin{array}{ccc}
3&7&1\\
7&21&0\\
1&0&3
\end{array}\!\right)$  &  $\mathcal C_{14} \cap \mathcal C_8 $  &  $1 = 54 - (45 + 17 - 9)$ & $45, 45$ & $18,19$ & $9$ & $1$ & $ d_1(\mu) = 3 $ &  $\begin{array}{c} W\subset\mathbb{P}^{8} \\ \deg(W)=6, g(W)=1 \\ \operatorname{gens\,id.}: {2}^{9}\  \end{array} $  &  $\begin{array}{c} \deg(U)=21, g(U)=16 \\ \chi({\mathcal O}_U)= 2 \\ \operatorname{gens\,id.}: {2}^{10}\ {3}^{8}\  \end{array} $  &  $\begin{array}{ll} \operatorname{num\,lines}: & 9 \\ \deg\operatorname{others}: & 0 \end{array}$  &  $\begin{array}{c} g(\widetilde{U}) = 16 \\ \deg(\widetilde{U})= 30 \end{array} $  &  $\det \left(\!\begin{array}{cc}
30&9\\
9&2
\end{array}\!\right) = -21 $  & 1\\
\hline
2 & $\mathfrak{S}\left(\!\begin{array}{cccc}
4&8&0&0\\
3&8&0&0
\end{array}\!\right) $ & $\begin{array}{c} \deg(S)=8 , g(S)= 3 \\ \operatorname{gens\,id.}: {2}^{2}\ {3}^{3}\ {4}^{1}\  \\ \deg(C)= 4 , n(C)= 2 \end{array}$  & $29 = 3 \cdot 8 + 5 =  \det \left(\!\begin{array}{ccc}
3&8&1\\
8&30&-1\\
1&-1&3
\end{array}\!\right)$  &  $\mathcal C_{26} \cap \mathcal C_8  \cap \mathcal C_{12} $  &  $1 = 54 - (47 + 14 - 8)$ & $47, 47$ & $15,15$ & $8$ & $2$ & $ d_1(\mu) = 3 $ &  $\begin{array}{c} W\subset\mathbb{P}^{7} \\ \deg(W)=5, g(W)=1 \\ \operatorname{gens\,id.}: {2}^{5}\  \end{array} $  &  $\begin{array}{c} \deg(U)=18, g(U)=14 \\ \chi({\mathcal O}_U)= 2 \\ \operatorname{gens\,id.}: {2}^{6}\ {3}^{6}\  \end{array} $  &  $\begin{array}{ll} \operatorname{num\,lines}: & 8 \\ \deg\operatorname{others}: & 0 \end{array}$  &  $\begin{array}{c} g(\widetilde{U}) = 14 \\ \deg(\widetilde{U})= 26 \end{array} $  &  $\det \left(\!\begin{array}{cc}
26&9\\
9&2
\end{array}\!\right) = -29 $  & 2\\
\hline
3 & $\mathfrak{S}\left(\!\begin{array}{cccc}
2&0&0&0\\
1&0&0&0
\end{array}\!\right) $ & $\begin{array}{c} \deg(S)=4 , g(S)= 0 \\ \operatorname{gens\,id.}: {2}^{6}\  \\ \deg(C)= 2 , n(C)= 0 \end{array}$  & $37 = 4 \cdot 8 + 5 =  \det \left(\!\begin{array}{ccc}
3&4&1\\
4&12&-1\\
1&-1&3
\end{array}\!\right)$  &  $\mathcal C_{20} \cap \mathcal C_8  \cap \mathcal C_{14} $  &  $1 = 54 - (29 + 24 - 0)$ & $29, 27$ & $25,28$ & $0$ & $1$ & $ d_1(\mu) = 2 $ &  $W=\mathbb{P}^4$  &  $\begin{array}{c} \deg(U)=6, g(U)=4 \\ \chi({\mathcal O}_U)= 2 \\ \operatorname{gens\,id.}: {2}^{1}\ {3}^{1}\  \end{array} $  &  $\begin{array}{ll} \operatorname{num\,lines}: & 0 \\ \deg\operatorname{others}: & 0 \end{array}$  &  $\begin{array}{c} g(\widetilde{U}) = 4 \\ \deg(\widetilde{U})= 6 \end{array} $  &  $\det \left(\!\begin{array}{cc}
6&7\\
7&2
\end{array}\!\right) = -37 $  & 3\\
\hline
4 & $\mathfrak{S}\left(\!\begin{array}{cccc}
3&3&0&0\\
1&0&0&0
\end{array}\!\right) $ & $\begin{array}{c} \deg(S)=6 , g(S)= 1 \\ \operatorname{gens\,id.}: {2}^{3}\ {3}^{2}\  \\ \deg(C)= 3 , n(C)= 1 \end{array}$  & $37 = 4 \cdot 8 + 5 =  \det \left(\!\begin{array}{ccc}
3&6&1\\
6&20&-1\\
1&-1&3
\end{array}\!\right)$  &  $\mathcal C_{24} \cap \mathcal C_8 $  &  $1 = 54 - (38 + 18 - 3)$ & $38, 38$ & $19,20$ & $3$ & $1$ & $ d_1(\mu) = 3 $ &  $\begin{array}{c} W\subset\mathbb{P}^{9} \\ \deg(W)=9, g(W)=3 \\ \operatorname{gens\,id.}: {2}^{12}\  \end{array} $  &  $\begin{array}{c} \deg(U)=19, g(U)=12 \\ \chi({\mathcal O}_U)= 2 \\ \operatorname{gens\,id.}: {2}^{18}\  \end{array} $  &  $\begin{array}{ll} \operatorname{num\,lines}: & 3 \\ \deg\operatorname{others}: & 0 \end{array}$  &  $\begin{array}{c} g(\widetilde{U}) = 12 \\ \deg(\widetilde{U})= 22 \end{array} $  &  $\det \left(\!\begin{array}{cc}
22&9\\
9&2
\end{array}\!\right) = -37 $  & 4\\
\hline
5 & $\mathfrak{S}\left(\!\begin{array}{cccc}
5&8&2&0\\
3&7&2&0
\end{array}\!\right) $ & $\begin{array}{c} \deg(S)=9 , g(S)= 4 \\ \operatorname{gens\,id.}: {2}^{1}\ {3}^{6}\ {4}^{1}\  \\ \deg(C)= 4 , n(C)= 2 \end{array}$  & $37 = 4 \cdot 8 + 5 =  \det \left(\!\begin{array}{ccc}
3&9&1\\
9&35&0\\
1&0&3
\end{array}\!\right)$  &  $\mathcal C_{24} \cap \mathcal C_8  \cap \mathcal C_{14} $  &  $1 = 54 - (50 + 11 - 8)$ & $50, 51$ & $12,12$ & $8$ & $1$ & $ d_1(\mu) = 3 $ &  $\begin{array}{c} W\subset\mathbb{P}^{5} \\ \deg(W)=2, g(W)=0 \\ \operatorname{gens\,id.}: {2}^{1}\  \end{array} $  &  $\begin{array}{c} \deg(U)=14, g(U)=12 \\ \chi({\mathcal O}_U)= -2 \\ \operatorname{gens\,id.}: {2}^{1}\ {3}^{1}\ {4}^{5}\ {5}^{1}\  \end{array} $  &  $\begin{array}{ll} \operatorname{num\,lines}: & 8 \\ \deg\operatorname{others}: & 0 \end{array}$  &  $\begin{array}{c} g(\widetilde{U}) = 12 \\ \deg(\widetilde{U})= 22 \end{array} $  &  $\det \left(\!\begin{array}{cc}
22&9\\
9&2
\end{array}\!\right) = -37 $  & 5\\
\hline
6 & $\mathfrak{S}\left(\!\begin{array}{cccc}
4&7&0&0\\
3&7&0&0
\end{array}\!\right) $ & $\begin{array}{c} \deg(S)=9 , g(S)= 3 \\ \operatorname{gens\,id.}: {2}^{1}\ {3}^{6}\ {5}^{1}\  \\ \deg(C)= 5 , n(C)= 5 \end{array}$  & $37 = 4 \cdot 8 + 5 =  \det \left(\!\begin{array}{ccc}
3&9&1\\
9&41&-2\\
1&-2&3
\end{array}\!\right)$  &  $\mathcal C_{42} \cap \mathcal C_8  \cap \mathcal C_{20} $  &  $1 = 54 - (49 + 11 - 7)$ & $49, 49$ & $12,12$ & $7$ & $2$ & $ d_1(\mu) = 3 $ &  $\begin{array}{c} W\subset\mathbb{P}^{6} \\ \deg(W)=4, g(W)=1 \\ \operatorname{gens\,id.}: {2}^{2}\  \end{array} $  &  $\begin{array}{c} \deg(U)=15, g(U)=12 \\ \chi({\mathcal O}_U)= 2 \\ \operatorname{gens\,id.}: {2}^{3}\ {3}^{4}\  \end{array} $  &  $\begin{array}{ll} \operatorname{num\,lines}: & 7 \\ \deg\operatorname{others}: & 0 \end{array}$  &  $\begin{array}{c} g(\widetilde{U}) = 12 \\ \deg(\widetilde{U})= 22 \end{array} $  &  $\det \left(\!\begin{array}{cc}
22&9\\
9&2
\end{array}\!\right) = -37 $  & 6\\
\hline
7 & $\mathfrak{S}\left(\!\begin{array}{cccc}
4&4&1&0\\
1&0&0&0
\end{array}\!\right) $ & $\begin{array}{c} \deg(S)=8 , g(S)= 2 \\ \operatorname{gens\,id.}: {2}^{1}\ {3}^{7}\ {4}^{1}\  \\ \deg(C)= 4 , n(C)= 3 \end{array}$  & $45 = 5 \cdot 8 + 5 =  \det \left(\!\begin{array}{ccc}
3&8&1\\
8&32&-1\\
1&-1&3
\end{array}\!\right)$  &  $\mathcal C_{32} \cap \mathcal C_8 $  &  $1 = 54 - (45 + 12 - 4)$ & $45, 46$ & $13,13$ & $4$ & $1$ & $ d_1(\mu) = 3 $ &  $\begin{array}{c} W\subset\mathbb{P}^{6} \\ \deg(W)=4, g(W)=1 \\ \operatorname{gens\,id.}: {2}^{2}\  \end{array} $  &  $\begin{array}{c} \deg(U)=14, g(U)=10 \\ \chi({\mathcal O}_U)= 1 \\ \operatorname{gens\,id.}: {2}^{3}\ {3}^{5}\  \end{array} $  &  $\begin{array}{ll} \operatorname{num\,lines}: & 4 \\ \deg\operatorname{others}: & 0 \end{array}$  &  $\begin{array}{c} g(\widetilde{U}) = 10 \\ \deg(\widetilde{U})= 18 \end{array} $  &  $\det \left(\!\begin{array}{cc}
18&9\\
9&2
\end{array}\!\right) = -45 $  & 7\\
\hline
8 & $\mathfrak{S}\left(\!\begin{array}{cccc}
5&10&1&0\\
3&8&1&0
\end{array}\!\right) $ & $\begin{array}{c} \deg(S)=11 , g(S)= 5 \\ \operatorname{gens\,id.}: {3}^{6}\ {4}^{2}\ {5}^{1}\  \\ \deg(C)= 5 , n(C)= 5 \end{array}$  & $45 = 5 \cdot 8 + 5 =  \det \left(\!\begin{array}{ccc}
3&11&1\\
11&51&0\\
1&0&3
\end{array}\!\right)$  &  $\mathcal C_{32} \cap \mathcal C_8  \cap \mathcal C_{18} $  &  $1 = 54 - (55 + 5 - 7)$ & $55,  - $ & $6,6$ & $7$ &  -  & $ d_1(\mu) = 4 $ &  $\begin{array}{c} W\subset\mathbb{P}^{9} \\ \deg(W)=9, g(W)=3 \\ \operatorname{gens\,id.}: {2}^{12}\  \end{array} $  &  $\begin{array}{c} \deg(U)=26, g(U)=19 \\ \chi({\mathcal O}_U)= 0 \\ \operatorname{gens\,id.}: {2}^{13}\ {3}^{13}\  \end{array} $  &  $\begin{array}{ll} \operatorname{num\,lines}: & 8 \\ \deg\operatorname{others}: & 2 \end{array}$  &  $\begin{array}{c} g(\widetilde{U}) = 20 \\ \deg(\widetilde{U})= 38 \end{array} $  &  $\det \left(\!\begin{array}{cc}
38&11\\
11&2
\end{array}\!\right) = -45 $  & 8\\
\hline
9 & $\mathfrak{S}\left(\!\begin{array}{cccc}
3&2&0&0\\
2&2&0&0
\end{array}\!\right) $ & $\begin{array}{c} \deg(S)=7 , g(S)= 1 \\ \operatorname{gens\,id.}: {2}^{2}\ {3}^{4}\ {4}^{1}\  \\ \deg(C)= 4 , n(C)= 3 \end{array}$  & $45 = 5 \cdot 8 + 5 =  \det \left(\!\begin{array}{ccc}
3&7&1\\
7&29&-2\\
1&-2&3
\end{array}\!\right)$  &  $\mathcal C_{38} \cap \mathcal C_8 $  &  $1 = 54 - (40 + 15 - 2)$ & $40, 40$ & $16,16$ & $2$ & $2$ & $ d_1(\mu) = 3 $ &  $\begin{array}{c} W\subset\mathbb{P}^{8} \\ \deg(W)=8, g(W)=3 \\ \operatorname{gens\,id.}: {2}^{7}\  \end{array} $  &  $\begin{array}{c} \deg(U)=16, g(U)=10 \\ \chi({\mathcal O}_U)= 2 \\ \operatorname{gens\,id.}: {2}^{13}\  \end{array} $  &  $\begin{array}{ll} \operatorname{num\,lines}: & 2 \\ \deg\operatorname{others}: & 0 \end{array}$  &  $\begin{array}{c} g(\widetilde{U}) = 10 \\ \deg(\widetilde{U})= 18 \end{array} $  &  $\det \left(\!\begin{array}{cc}
18&9\\
9&2
\end{array}\!\right) = -45 $  & 9\\
\hline
10 & $\mathfrak{S}\left(\!\begin{array}{cccc}
5&6&0&1\\
1&0&0&0
\end{array}\!\right) $ & $\begin{array}{c} \deg(S)=10 , g(S)= 3 \\ \operatorname{gens\,id.}: {3}^{7}\ {4}^{1}\ {5}^{1}\  \\ \deg(C)= 5 , n(C)= 6 \end{array}$  & $45 = 5 \cdot 8 + 5 =  \det \left(\!\begin{array}{ccc}
3&10&1\\
10&46&-1\\
1&-1&3
\end{array}\!\right)$  &  $\mathcal C_{38} \cap \mathcal C_8 $  &  $1 = 54 - (52 + 6 - 5)$ & $52,  - $ & $7,7$ & $5$ &  -  & $ d_1(\mu) = 4 $ &  $\begin{array}{c} W\subset\mathbb{P}^{11} \\ \deg(W)=14, g(W)=6 \\ \operatorname{gens\,id.}: {2}^{22}\  \end{array} $  &  $\begin{array}{c} \deg(U)=28, g(U)=19 \\ \chi({\mathcal O}_U)= 2 \\ \operatorname{gens\,id.}: {2}^{28}\ {3}^{2}\  \end{array} $  &  $\begin{array}{ll} \operatorname{num\,lines}: & 6 \\ \deg\operatorname{others}: & 2 \end{array}$  &  $\begin{array}{c} g(\widetilde{U}) = 20 \\ \deg(\widetilde{U})= 38 \end{array} $  &  $\det \left(\!\begin{array}{cc}
38&11\\
11&2
\end{array}\!\right) = -45 $  & 10\\
\hline
11 & $\mathfrak{S}\left(\!\begin{array}{cccc}
5&7&2&0\\
3&6&2&0
\end{array}\!\right) $ & $\begin{array}{c} \deg(S)=10 , g(S)= 4 \\ \operatorname{gens\,id.}: {3}^{9}\ {5}^{1}\  \\ \deg(C)= 5 , n(C)= 5 \end{array}$  & $45 = 5 \cdot 8 + 5 =  \det \left(\!\begin{array}{ccc}
3&10&1\\
10&46&-1\\
1&-1&3
\end{array}\!\right)$  &  $\mathcal C_{38} \cap \mathcal C_8  \cap \mathcal C_{20} $  &  $1 = 54 - (52 + 8 - 7)$ & $52, 52$ & $9,9$ & $7$ & $2$ & $ d_1(\mu) = 3 $ &  $W=\mathbb{P}^4$  &  $\begin{array}{c} \deg(U)=11, g(U)=10 \\ \chi({\mathcal O}_U)= -5 \\ \operatorname{gens\,id.}: {4}^{1}\ {5}^{6}\  \end{array} $  &  $\begin{array}{ll} \operatorname{num\,lines}: & 7 \\ \deg\operatorname{others}: & 0 \end{array}$  &  $\begin{array}{c} g(\widetilde{U}) = 10 \\ \deg(\widetilde{U})= 18 \end{array} $  &  $\det \left(\!\begin{array}{cc}
18&9\\
9&2
\end{array}\!\right) = -45 $  & 11\\
\hline
12 & $\mathfrak{S}\left(\!\begin{array}{cccc}
4&6&0&0\\
3&6&0&0
\end{array}\!\right) $ & $\begin{array}{c} \deg(S)=10 , g(S)= 3 \\ \operatorname{gens\,id.}: {3}^{9}\ {6}^{1}\  \\ \deg(C)= 6 , n(C)= 9 \end{array}$  & $45 = 5 \cdot 8 + 5 =  \det \left(\!\begin{array}{ccc}
3&10&1\\
10&54&-3\\
1&-3&3
\end{array}\!\right)$  &  $\mathcal C_{62} \cap \mathcal C_8  \cap \mathcal C_{32} $  &  $1 = 54 - (51 + 8 - 6)$ & $51, 51$ & $9,9$ & $6$ & $2$(!)  & $ d_1(\mu) = 3 $ &  $\begin{array}{c} W\subset\mathbb{P}^{5} \\ \deg(W)=3, g(W)=1 \\ \operatorname{gens\,id.}: {3}^{1}\  \end{array} $  &  $\begin{array}{c} \deg(U)=12, g(U)=10 \\ \chi({\mathcal O}_U)= 2 \\ \operatorname{gens\,id.}: {2}^{1}\ {3}^{3}\  \end{array} $  &  $\begin{array}{ll} \operatorname{num\,lines}: & 6 \\ \deg\operatorname{others}: & 0 \end{array}$  &  $\begin{array}{c} g(\widetilde{U}) = 10 \\ \deg(\widetilde{U})= 18 \end{array} $  &  $\det \left(\!\begin{array}{cc}
18&9\\
9&2
\end{array}\!\right) = -45 $  & 12\\
\hline
13 & $\mathfrak{S}\left(\!\begin{array}{cccc}
4&3&1&0\\
2&1&1&0
\end{array}\!\right) $ & $\begin{array}{c} \deg(S)=9 , g(S)= 2 \\ \operatorname{gens\,id.}: {3}^{10}\ {5}^{1}\  \\ \deg(C)= 5 , n(C)= 6 \end{array}$  & $53 = 6 \cdot 8 + 5 =  \det \left(\!\begin{array}{ccc}
3&9&1\\
9&43&-2\\
1&-2&3
\end{array}\!\right)$  &  $\mathcal C_{48} \cap \mathcal C_8 $  &  $1 = 54 - (47 + 9 - 3)$ & $47, 47$ & $10,10$ & $3$ & $2$ & $ d_1(\mu) = 3 $ &  $\begin{array}{c} W\subset\mathbb{P}^{5} \\ \deg(W)=3, g(W)=1 \\ \operatorname{gens\,id.}: {3}^{1}\  \end{array} $  &  $\begin{array}{c} \deg(U)=11, g(U)=8 \\ \chi({\mathcal O}_U)= 1 \\ \operatorname{gens\,id.}: {2}^{1}\ {3}^{4}\  \end{array} $  &  $\begin{array}{ll} \operatorname{num\,lines}: & 3 \\ \deg\operatorname{others}: & 0 \end{array}$  &  $\begin{array}{c} g(\widetilde{U}) = 8 \\ \deg(\widetilde{U})= 14 \end{array} $  &  $\det \left(\!\begin{array}{cc}
14&9\\
9&2
\end{array}\!\right) = -53 $  & 13\\
\hline
14 & $\mathfrak{S}\left(\!\begin{array}{cccc}
5&9&1&0\\
3&7&1&0
\end{array}\!\right) $ & $\begin{array}{c} \deg(S)=12 , g(S)= 5 \\ \operatorname{gens\,id.}: {3}^{3}\ {4}^{12}\ {6}^{1}\  \\ \deg(C)= 6 , n(C)= 9 \end{array}$  & $53 = 6 \cdot 8 + 5 =  \det \left(\!\begin{array}{ccc}
3&12&1\\
12&64&-1\\
1&-1&3
\end{array}\!\right)$  &  $\mathcal C_{48} \cap \mathcal C_8 $  &  $1 = 54 - (57 + 2 - 6)$ & $57,  - $ & $3,3$ & $6$ &  -  & $ d_1(\mu) = 4 $ &  $\begin{array}{c} W\subset\mathbb{P}^{8} \\ \deg(W)=8, g(W)=3 \\ \operatorname{gens\,id.}: {2}^{7}\  \end{array} $  &  $\begin{array}{c} \deg(U)=23, g(U)=17 \\ \chi({\mathcal O}_U)= 0 \\ \operatorname{gens\,id.}: {2}^{8}\ {3}^{11}\  \end{array} $  &  $\begin{array}{ll} \operatorname{num\,lines}: & 7 \\ \deg\operatorname{others}: & 2 \end{array}$  &  $\begin{array}{c} g(\widetilde{U}) = 18 \\ \deg(\widetilde{U})= 34 \end{array} $  &  $\det \left(\!\begin{array}{cc}
34&23\\
23&14
\end{array}\!\right) = -53 $  & 14\\
\hline
15 & $\mathfrak{S}\left(\!\begin{array}{cccc}
3&1&0&0\\
2&1&0&0
\end{array}\!\right) $ & $\begin{array}{c} \deg(S)=8 , g(S)= 1 \\ \operatorname{gens\,id.}: {2}^{1}\ {3}^{7}\ {5}^{1}\  \\ \deg(C)= 5 , n(C)= 6 \end{array}$  & $53 = 6 \cdot 8 + 5 =  \det \left(\!\begin{array}{ccc}
3&8&1\\
8&40&-3\\
1&-3&3
\end{array}\!\right)$  &  $\mathcal C_{56} \cap \mathcal C_8 $  &  $1 = 54 - (42 + 12 - 1)$ & $42, 42$ & $13,13$ & $1$ & $2$ & $ d_1(\mu) = 3 $ &  $\begin{array}{c} W\subset\mathbb{P}^{7} \\ \deg(W)=7, g(W)=3 \\ \operatorname{gens\,id.}: {2}^{3}\ {3}^{1}\  \end{array} $  &  $\begin{array}{c} \deg(U)=13, g(U)=8 \\ \chi({\mathcal O}_U)= 2 \\ \operatorname{gens\,id.}: {2}^{9}\  \end{array} $  &  $\begin{array}{ll} \operatorname{num\,lines}: & 1 \\ \deg\operatorname{others}: & 0 \end{array}$  &  $\begin{array}{c} g(\widetilde{U}) = 8 \\ \deg(\widetilde{U})= 14 \end{array} $  &  $\det \left(\!\begin{array}{cc}
14&9\\
9&2
\end{array}\!\right) = -53 $  & 15\\
\hline
16 & $\mathfrak{S}\left(\!\begin{array}{cccc}
5&5&0&1\\
2&1&0&1
\end{array}\!\right) $ & $\begin{array}{c} \deg(S)=11 , g(S)= 3 \\ \operatorname{gens\,id.}: {3}^{4}\ {4}^{9}\ {6}^{1}\  \\ \deg(C)= 6 , n(C)= 10 \end{array}$  & $53 = 6 \cdot 8 + 5 =  \det \left(\!\begin{array}{ccc}
3&11&1\\
11&59&-2\\
1&-2&3
\end{array}\!\right)$  &  $\mathcal C_{56} \cap \mathcal C_8 $  &  $1 = 54 - (54 + 3 - 4)$ & $54,  - $ & $4,4$ & $4$ &  -  & $ d_1(\mu) = 4 $ &  $\begin{array}{c} W\subset\mathbb{P}^{10} \\ \deg(W)=13, g(W)=6 \\ \operatorname{gens\,id.}: {2}^{15}\  \end{array} $  &  $\begin{array}{c} \deg(U)=25, g(U)=17 \\ \chi({\mathcal O}_U)= 2 \\ \operatorname{gens\,id.}: {2}^{21}\ {3}^{2}\  \end{array} $  &  $\begin{array}{ll} \operatorname{num\,lines}: & 5 \\ \deg\operatorname{others}: & 2 \end{array}$  &  $\begin{array}{c} g(\widetilde{U}) = 18 \\ \deg(\widetilde{U})= 34 \end{array} $  &  $\det \left(\!\begin{array}{cc}
34&11\\
11&2
\end{array}\!\right) = -53 $  & 16\\
\hline
17 & $\mathfrak{S}\left(\!\begin{array}{cccc}
5&6&2&0\\
3&5&2&0
\end{array}\!\right) $ & $\begin{array}{c} \deg(S)=11 , g(S)= 4 \\ \operatorname{gens\,id.}: {3}^{6}\ {4}^{2}\ {6}^{1}\  \\ \deg(C)= 6 , n(C)= 9 \end{array}$  & $53 = 6 \cdot 8 + 5 =  \det \left(\!\begin{array}{ccc}
3&11&1\\
11&59&-2\\
1&-2&3
\end{array}\!\right)$  &  $\mathcal C_{56} \cap \mathcal C_8  \cap \mathcal C_{30} $  &  $1 = 54 - (54 + 5 - 6)$ & $54,  - $ & $6,6$ & $6$ &  -  & $ d_1(\mu) = 4 $ &  $\begin{array}{c} W\subset\mathbb{P}^{11} \\ \deg(W)=14, g(W)=6 \\ \operatorname{gens\,id.}: {2}^{22}\  \end{array} $  &  $\begin{array}{c} \deg(U)=28, g(U)=18 \\ \chi({\mathcal O}_U)= 0 \\ \operatorname{gens\,id.}: {2}^{28}\  \end{array} $  &  $\begin{array}{ll} \operatorname{num\,lines}: & 6 \\ \deg\operatorname{others}: & 0 \end{array}$  &  $\begin{array}{c} g(\widetilde{U}) = 18 \\ \deg(\widetilde{U})= 34 \end{array} $  &  $\det \left(\!\begin{array}{cc}
34&11\\
11&2
\end{array}\!\right) = -53 $  & 17\\
\hline
18 & $\mathfrak{S}\left(\!\begin{array}{cccc}
4&5&0&0\\
3&5&0&0
\end{array}\!\right) $ & $\begin{array}{c} \deg(S)=11 , g(S)= 3 \\ \operatorname{gens\,id.}: {3}^{6}\ {4}^{3}\ {7}^{1}\  \\ \deg(C)= 7 , n(C)= 14 \end{array}$  & $53 = 6 \cdot 8 + 5 =  \det \left(\!\begin{array}{ccc}
3&11&1\\
11&69&-4\\
1&-4&3
\end{array}\!\right)$  &  $\mathcal C_{86} \cap \mathcal C_8 $  &  $1 = 54 - (53 + 5 - 5)$ & $53,  - $ & $6,6$ & $5$ &  -  & $ d_1(\mu) = 4 $ &  $\begin{array}{c} W\subset\mathbb{P}^{13} \\ \deg(W)=20, g(W)=10 \\ \operatorname{gens\,id.}: {2}^{35}\  \end{array} $  &  $\begin{array}{c} \deg(U)=29, g(U)=18 \\ \chi({\mathcal O}_U)= 2 \\ \operatorname{gens\,id.}: {2}^{50}\  \end{array} $  &  $\begin{array}{ll} \operatorname{num\,lines}: & 5 \\ \deg\operatorname{others}: & 0 \end{array}$  &  $\begin{array}{c} g(\widetilde{U}) = 18 \\ \deg(\widetilde{U})= 34 \end{array} $  &  $\det \left(\!\begin{array}{cc}
34&11\\
11&2
\end{array}\!\right) = -53 $  & 18\\
\hline
19 & $\mathfrak{S}\left(\!\begin{array}{cccc}
6&5&5&0\\
3&5&4&0
\end{array}\!\right) $ & $\begin{array}{c} \deg(S)=11 , g(S)= 5 \\ \operatorname{gens\,id.}: {3}^{6}\ {4}^{2}\ {5}^{1}\  \\ \deg(C)= 5 , n(C)= 5 \end{array}$  & $61 = 7 \cdot 8 + 5 =  \det \left(\!\begin{array}{ccc}
3&11&1\\
11&53&0\\
1&0&3
\end{array}\!\right)$  &  $\mathcal C_{38} \cap \mathcal C_8  \cap \mathcal C_{24} $  &  $1 = 54 - (53 + 5 - 5)$ & $53, 54$ & $6,6$ & $5$ & $1$ & $ d_1(\mu) = 4 $ &  $\begin{array}{c} W\subset\mathbb{P}^{9} \\ \deg(W)=9, g(W)=3 \\ \operatorname{gens\,id.}: {2}^{12}\  \end{array} $  &  $\begin{array}{c} \deg(U)=25, g(U)=16 \\ \chi({\mathcal O}_U)= -3 \\ \operatorname{gens\,id.}: {2}^{13}\ {3}^{13}\  \end{array} $  &  $\begin{array}{ll} \operatorname{num\,lines}: & 5 \\ \deg\operatorname{others}: & 0 \end{array}$  &  $\begin{array}{c} g(\widetilde{U}) = 16 \\ \deg(\widetilde{U})= 30 \end{array} $  &  $\det \left(\!\begin{array}{cc}
30&11\\
11&2
\end{array}\!\right) = -61 $  & 19\\
\hline
20 & $\mathfrak{S}\left(\!\begin{array}{cccc}
4&6&0&0\\
2&3&0&0
\end{array}\!\right) $ & $\begin{array}{c} \deg(S)=10 , g(S)= 3 \\ \operatorname{gens\,id.}: {3}^{7}\ {4}^{1}\ {5}^{1}\  \\ \deg(C)= 5 , n(C)= 6 \end{array}$  & $61 = 7 \cdot 8 + 5 =  \det \left(\!\begin{array}{ccc}
3&10&1\\
10&48&-1\\
1&-1&3
\end{array}\!\right)$  &  $\mathcal C_{44} \cap \mathcal C_8 $  &  $1 = 54 - (50 + 6 - 3)$ & $50, 51$ & $7,7$ & $3$ & $1$ & $ d_1(\mu) = 4 $ &  $\begin{array}{c} W\subset\mathbb{P}^{11} \\ \deg(W)=14, g(W)=6 \\ \operatorname{gens\,id.}: {2}^{22}\  \end{array} $  &  $\begin{array}{c} \deg(U)=27, g(U)=16 \\ \chi({\mathcal O}_U)= -1 \\ \operatorname{gens\,id.}: {2}^{28}\  \end{array} $  &  $\begin{array}{ll} \operatorname{num\,lines}: & 3 \\ \deg\operatorname{others}: & 0 \end{array}$  &  $\begin{array}{c} g(\widetilde{U}) = 16 \\ \deg(\widetilde{U})= 30 \end{array} $  &  $\det \left(\!\begin{array}{cc}
30&11\\
11&2
\end{array}\!\right) = -61 $  & 20\\
\hline
21 & $\mathfrak{S}\left(\!\begin{array}{cccc}
4&2&1&0\\
2&0&1&0
\end{array}\!\right) $ & $\begin{array}{c} \deg(S)=10 , g(S)= 2 \\ \operatorname{gens\,id.}: {3}^{7}\ {4}^{2}\ {6}^{1}\  \\ \deg(C)= 6 , n(C)= 10 \end{array}$  & $61 = 7 \cdot 8 + 5 =  \det \left(\!\begin{array}{ccc}
3&10&1\\
10&56&-3\\
1&-3&3
\end{array}\!\right)$  &  $\mathcal C_{68} \cap \mathcal C_8 $  &  $1 = 54 - (49 + 6 - 2)$ & $49, 49$ & $7,7$ & $2$ & $2$(!)  & $ d_1(\mu) = 4 $ &  $\begin{array}{c} W\subset\mathbb{P}^{13} \\ \deg(W)=20, g(W)=10 \\ \operatorname{gens\,id.}: {2}^{35}\  \end{array} $  &  $\begin{array}{c} \deg(U)=28, g(U)=16 \\ \chi({\mathcal O}_U)= 1 \\ \operatorname{gens\,id.}: {2}^{50}\  \end{array} $  &  $\begin{array}{ll} \operatorname{num\,lines}: & 2 \\ \deg\operatorname{others}: & 0 \end{array}$  &  $\begin{array}{c} g(\widetilde{U}) = 16 \\ \deg(\widetilde{U})= 30 \end{array} $  &  $\det \left(\!\begin{array}{cc}
30&11\\
11&2
\end{array}\!\right) = -61 $  & 21\\
\hline
22 & $\mathfrak{S}\left(\!\begin{array}{cccc}
3&0&0&0\\
2&0&0&0
\end{array}\!\right) $ & $\begin{array}{c} \deg(S)=9 , g(S)= 1 \\ \operatorname{gens\,id.}: {3}^{10}\ {6}^{1}\  \\ \deg(C)= 6 , n(C)= 10 \end{array}$  & $61 = 7 \cdot 8 + 5 =  \det \left(\!\begin{array}{ccc}
3&9&1\\
9&53&-4\\
1&-4&3
\end{array}\!\right)$  &  $\mathcal C_{78} \cap \mathcal C_8 $  &  $1 = 54 - (44 + 9 - 0)$ & $44, 44$ & $10,10$ & $0$ & $2$(!)  & $ d_1(\mu) = 3 $ &  $\begin{array}{c} W\subset\mathbb{P}^{6} \\ \deg(W)=6, g(W)=3 \\ \operatorname{gens\,id.}: {3}^{4}\  \end{array} $  &  $\begin{array}{c} \deg(U)=10, g(U)=6 \\ \chi({\mathcal O}_U)= 2 \\ \operatorname{gens\,id.}: {2}^{6}\  \end{array} $  &  $\begin{array}{ll} \operatorname{num\,lines}: & 0 \\ \deg\operatorname{others}: & 0 \end{array}$  &  $\begin{array}{c} g(\widetilde{U}) = 6 \\ \deg(\widetilde{U})= 10 \end{array} $  &  $\det \left(\!\begin{array}{cc}
10&9\\
9&2
\end{array}\!\right) = -61 $  & 22\\
\hline
23 & $\mathfrak{S}\left(\!\begin{array}{cccc}
6&4&5&0\\
3&4&4&0
\end{array}\!\right) $ & $\begin{array}{c} \deg(S)=12 , g(S)= 5 \\ \operatorname{gens\,id.}: {3}^{3}\ {4}^{12}\ {6}^{1}\  \\ \deg(C)= 6 , n(C)= 9 \end{array}$  & $69 = 8 \cdot 8 + 5 =  \det \left(\!\begin{array}{ccc}
3&12&1\\
12&66&-1\\
1&-1&3
\end{array}\!\right)$  &  $\mathcal C_{54} \cap \mathcal C_8 $  &  $1 = 54 - (55 + 2 - 4)$ & $55, 55$ & $3,3$ & $4$ & $2$(!)  & $ d_1(\mu) = 4 $ &  $\begin{array}{c} W\subset\mathbb{P}^{8} \\ \deg(W)=8, g(W)=3 \\ \operatorname{gens\,id.}: {2}^{7}\  \end{array} $  &  $\begin{array}{c} \deg(U)=22, g(U)=14 \\ \chi({\mathcal O}_U)= -3 \\ \operatorname{gens\,id.}: {2}^{8}\ {3}^{11}\  \end{array} $  &  $\begin{array}{ll} \operatorname{num\,lines}: & 4 \\ \deg\operatorname{others}: & 0 \end{array}$  &  $\begin{array}{c} g(\widetilde{U}) = 14 \\ \deg(\widetilde{U})= 26 \end{array} $  &  $\det \left(\!\begin{array}{cc}
26&11\\
11&2
\end{array}\!\right) = -69 $  & 23\\
\hline
24 & $\mathfrak{S}\left(\!\begin{array}{cccc}
4&5&0&0\\
2&2&0&0
\end{array}\!\right) $ & $\begin{array}{c} \deg(S)=11 , g(S)= 3 \\ \operatorname{gens\,id.}: {3}^{4}\ {4}^{9}\ {6}^{1}\  \\ \deg(C)= 6 , n(C)= 10 \end{array}$  & $69 = 8 \cdot 8 + 5 =  \det \left(\!\begin{array}{ccc}
3&11&1\\
11&61&-2\\
1&-2&3
\end{array}\!\right)$  &  $\mathcal C_{62} \cap \mathcal C_8 $  &  $1 = 54 - (52 + 3 - 2)$ & $52, 52$ & $4,4$ & $2$ & $2$(!)  & $ d_1(\mu) = 4 $ &  $\begin{array}{c} W\subset\mathbb{P}^{10} \\ \deg(W)=13, g(W)=6 \\ \operatorname{gens\,id.}: {2}^{15}\  \end{array} $  &  $\begin{array}{c} \deg(U)=24, g(U)=14 \\ \chi({\mathcal O}_U)= -1 \\ \operatorname{gens\,id.}: {2}^{21}\  \end{array} $  &  $\begin{array}{ll} \operatorname{num\,lines}: & 2 \\ \deg\operatorname{others}: & 0 \end{array}$  &  $\begin{array}{c} g(\widetilde{U}) = 14 \\ \deg(\widetilde{U})= 26 \end{array} $  &  $\det \left(\!\begin{array}{cc}
26&11\\
11&2
\end{array}\!\right) = -69 $  & 24
\\ \hline
\end{tabular}

\end{adjustbox}}
\restoregeometry
\clearpage
\pagestyle{headings}

Table \ref{TABex2} presents some additional examples obtained from the same constructions. For these examples, the parameter count gives an upper bound equal to \(2\) or \(3\) for \(\codim_{\mathcal C_8}\mathcal C_{\mathfrak S}\); see, however, Remark \ref{refinedCount}. As for some examples in Table \ref{TABex}, the invariants of \(\widetilde U\) and its lattice polarization are obtained virtually in Examples 2 and 5--12. In the last example, the final computations related to the associated K3 surface were not completed because of their computational cost.

For the current versions of the article and the archive, the identifiers corresponding to the rows of Table \ref{TABex2} are
\[ \begin{array}{c|rrrrrrrrrrrrr} i&1&2&3&4&5&6&7&8&9&10&11&12&13\\ \mathrm{ID}&2&1&4&7&10&8&15&13&24&18&20&29&27 \end{array}. \]

\newgeometry{
  left=0.5cm,
  right=0.5cm,
  top=0.5cm,
  bottom=0.5cm
}
\pagestyle{empty}
\thispagestyle{empty}
\refstepcounter{table}
\label{TABex2}
\vspace*{-\topskip}
\noindent
\rotatebox{90}{%
\begin{adjustbox}{
  width=0.98\textheight,
  totalheight=0.98\textwidth,
  keepaspectratio
}
\begin{tabular}{|c||c|c|l|c||l|c|c|c|c||c|c|c|c|c|r||c|}
\hline
$i$ & surface $S\subset X\subset \mathbb{P}^5$ & $S, C=S\cap P$ & & & parameter count &(\textdagger) & (\textdaggerdbl) & (\textsection) & (\textasteriskcentered) & $\mu:\mathbb{P}^5\dashrightarrow W$ & $W$ & $U$ & exc. curves & $\widetilde{U}$ & & $i$ \\
\hline
\\
\hline
1 & $\mathfrak{S}\left(\!\begin{array}{cccc}
3&4&0&0\\
1&1&0&0
\end{array}\!\right) $ & $\begin{array}{c} \deg(S)=5 , g(S)= 1 \\ \operatorname{gens\,id.}: {2}^{5}\  \\ \deg(C)= 2 , n(C)= 0 \end{array}$  & $29 = 3 \cdot 8 + 5 =  \det \left(\!\begin{array}{ccc}
3&5&1\\
5&13&0\\
1&0&3
\end{array}\!\right)$  &  $\mathcal C_{14} \cap \mathcal C_8 $  &  $2 = 54 - (36 + 21 - 5)$ & $36, 35$ & $22,25$ & $5$ & $1$ & $ d_1(\mu) = 3 $ &  $\begin{array}{c} W\subset\mathbb{P}^{10} \\ \deg(W)=10, g(W)=3 \\ \operatorname{gens\,id.}: {2}^{18}\  \end{array} $  &  $\begin{array}{c} \deg(U)=22, g(U)=14 \\ \chi({\mathcal O}_U)= 2 \\ \operatorname{gens\,id.}: {2}^{24}\  \end{array} $  &  $\begin{array}{ll} \operatorname{num\,lines}: & 4 \\ \deg\operatorname{others}: & 0 \end{array}$  &  $\begin{array}{c} g(\widetilde{U}) = 14 \\ \deg(\widetilde{U})= 26 \end{array} $  &  $\det \left(\!\begin{array}{cc}
26&9\\
9&2
\end{array}\!\right) = -29 $ & 1\\
\hline
2 & $\mathfrak{S}\left(\!\begin{array}{cccc}
4&6&1&0\\
1&2&0&0
\end{array}\!\right) $ & $\begin{array}{c} \deg(S)=6 , g(S)= 2 \\ \operatorname{gens\,id.}: {2}^{4}\  \\ \deg(C)= 2 , n(C)= 0 \end{array}$  & $29 = 3 \cdot 8 + 5 =  \det \left(\!\begin{array}{ccc}
3&6&1\\
6&16&1\\
1&1&3
\end{array}\!\right)$  &  $\mathcal C_{12} \cap \mathcal C_8 $  &  $3 = 54 - (41 + 18 - 8)$ & $41, 41$ & $19,22$ & $8$ & $1$ & $ d_1(\mu) = 3 $ &  $\begin{array}{c} W\subset\mathbb{P}^{8} \\ \deg(W)=6, g(W)=1 \\ \operatorname{gens\,id.}: {2}^{9}\  \end{array} $  &  $\begin{array}{c} \deg(U)=20, g(U)=14 \\ \chi({\mathcal O}_U)= 1 \\ \operatorname{gens\,id.}: {2}^{10}\ {3}^{9}\  \end{array} $  &  $\begin{array}{ll} \operatorname{num\,lines}: & 6 \\ \deg\operatorname{others}: & 0 \end{array}$  &  $\begin{array}{c} g(\widetilde{U}) = 14 \\ \deg(\widetilde{U})= 26 \end{array} $  &  $\det \left(\!\begin{array}{cc}
26&9\\
9&2
\end{array}\!\right) = -29 $ & 2\\
\hline
3 & $\mathfrak{S}\left(\!\begin{array}{cccc}
5&8&0&1\\
1&2&0&0
\end{array}\!\right) $ & $\begin{array}{c} \deg(S)=8 , g(S)= 3 \\ \operatorname{gens\,id.}: {2}^{1}\ {3}^{8}\  \\ \deg(C)= 3 , n(C)= 1 \end{array}$  & $29 = 3 \cdot 8 + 5 =  \det \left(\!\begin{array}{ccc}
3&8&1\\
8&26&1\\
1&1&3
\end{array}\!\right)$  &  $\mathcal C_{14} \cap \mathcal C_8 $  &  $3 = 54 - (48 + 12 - 9)$ & $48, 50$ & $13,14$ & $9$ & $1$ & $ d_1(\mu) = 3 $ &  $\begin{array}{c} W\subset\mathbb{P}^{5} \\ \deg(W)=2, g(W)=0 \\ \operatorname{gens\,id.}: {2}^{1}\  \end{array} $  &  $\begin{array}{c} \deg(U)=14, g(U)=13 \\ \chi({\mathcal O}_U)= 1 \\ \operatorname{gens\,id.}: {2}^{1}\ {3}^{1}\ {4}^{6}\ {5}^{1}\  \end{array} $  &  $\begin{array}{ll} \operatorname{num\,lines}: & 8 \\ \deg\operatorname{others}: & 2 \end{array}$  &  $\begin{array}{c} g(\widetilde{U}) = 14 \\ \deg(\widetilde{U})= 26 \end{array} $  &  $\det \left(\!\begin{array}{cc}
26&9\\
9&2
\end{array}\!\right) = -29 $ & 3\\
\hline
4 & $\mathfrak{S}\left(\!\begin{array}{cccc}
4&5&1&0\\
1&1&0&0
\end{array}\!\right) $ & $\begin{array}{c} \deg(S)=7 , g(S)= 2 \\ \operatorname{gens\,id.}: {2}^{2}\ {3}^{5}\  \\ \deg(C)= 3 , n(C)= 1 \end{array}$  & $37 = 4 \cdot 8 + 5 =  \det \left(\!\begin{array}{ccc}
3&7&1\\
7&23&0\\
1&0&3
\end{array}\!\right)$  &  $\mathcal C_{20} \cap \mathcal C_8 $  &  $2 = 54 - (43 + 15 - 6)$ & $43, 44$ & $16,17$ & $6$ & $1$ & $ d_1(\mu) = 3 $ &  $\begin{array}{c} W\subset\mathbb{P}^{7} \\ \deg(W)=5, g(W)=1 \\ \operatorname{gens\,id.}: {2}^{5}\  \end{array} $  &  $\begin{array}{c} \deg(U)=17, g(U)=12 \\ \chi({\mathcal O}_U)= 1 \\ \operatorname{gens\,id.}: {2}^{6}\ {3}^{7}\  \end{array} $  &  $\begin{array}{ll} \operatorname{num\,lines}: & 5 \\ \deg\operatorname{others}: & 0 \end{array}$  &  $\begin{array}{c} g(\widetilde{U}) = 12 \\ \deg(\widetilde{U})= 22 \end{array} $  &  $\det \left(\!\begin{array}{cc}
22&9\\
9&2
\end{array}\!\right) = -37 $ & 4\\
\hline
5 & $\mathfrak{S}\left(\!\begin{array}{cccc}
5&7&0&1\\
1&1&0&0
\end{array}\!\right) $ & $\begin{array}{c} \deg(S)=9 , g(S)= 3 \\ \operatorname{gens\,id.}: {3}^{10}\ {4}^{1}\  \\ \deg(C)= 4 , n(C)= 3 \end{array}$  & $37 = 4 \cdot 8 + 5 =  \det \left(\!\begin{array}{ccc}
3&9&1\\
9&35&0\\
1&0&3
\end{array}\!\right)$  &  $\mathcal C_{24} \cap \mathcal C_8 $  &  $2 = 54 - (50 + 9 - 7)$ & $50, 52$ & $10,10$ & $7$ & $1$ & $ d_1(\mu) = 3 $ &  $W=\mathbb{P}^4$  &  $\begin{array}{c} \deg(U)=11, g(U)=11 \\ \chi({\mathcal O}_U)= -1 \\ \operatorname{gens\,id.}: {4}^{1}\ {5}^{7}\  \end{array} $  &  $\begin{array}{ll} \operatorname{num\,lines}: & 7 \\ \deg\operatorname{others}: & 2 \end{array}$  &  $\begin{array}{c} g(\widetilde{U}) = 12 \\ \deg(\widetilde{U})= 22 \end{array} $  &  $\det \left(\!\begin{array}{cc}
22&9\\
9&2
\end{array}\!\right) = -37 $ & 5\\
\hline
6 & $\mathfrak{S}\left(\!\begin{array}{cccc}
4&8&0&0\\
2&5&0&0
\end{array}\!\right) $ & $\begin{array}{c} \deg(S)=8 , g(S)= 3 \\ \operatorname{gens\,id.}: {2}^{1}\ {3}^{8}\  \\ \deg(C)= 3 , n(C)= 1 \end{array}$  & $45 = 5 \cdot 8 + 5 =  \det \left(\!\begin{array}{ccc}
3&8&1\\
8&28&1\\
1&1&3
\end{array}\!\right)$  &  $\mathcal C_{20} \cap \mathcal C_8 $  &  $3 = 54 - (46 + 12 - 7)$ & $46, 48$ & $13,14$ & $7$ & $1$ & $ d_1(\mu) = 3 $ &  $\begin{array}{c} W\subset\mathbb{P}^{5} \\ \deg(W)=2, g(W)=0 \\ \operatorname{gens\,id.}: {2}^{1}\  \end{array} $  &  $\begin{array}{c} \deg(U)=13, g(U)=10 \\ \chi({\mathcal O}_U)= -2 \\ \operatorname{gens\,id.}: {2}^{1}\ {3}^{1}\ {4}^{7}\  \end{array} $  &  $\begin{array}{ll} \operatorname{num\,lines}: & 5 \\ \deg\operatorname{others}: & 0 \end{array}$  &  $\begin{array}{c} g(\widetilde{U}) = 10 \\ \deg(\widetilde{U})= 18 \end{array} $  &  $\det \left(\!\begin{array}{cc}
18&9\\
9&2
\end{array}\!\right) = -45 $ & 6\\
\hline
7 & $\mathfrak{S}\left(\!\begin{array}{cccc}
4&7&0&0\\
2&4&0&0
\end{array}\!\right) $ & $\begin{array}{c} \deg(S)=9 , g(S)= 3 \\ \operatorname{gens\,id.}: {3}^{10}\ {4}^{1}\  \\ \deg(C)= 4 , n(C)= 3 \end{array}$  & $53 = 6 \cdot 8 + 5 =  \det \left(\!\begin{array}{ccc}
3&9&1\\
9&37&0\\
1&0&3
\end{array}\!\right)$  &  $\mathcal C_{30} \cap \mathcal C_8 $  &  $2 = 54 - (48 + 9 - 5)$ & $48, 50$ & $10,10$ & $5$ & $1$ & $ d_1(\mu) = 3 $ &  $W=\mathbb{P}^4$  &  $\begin{array}{c} \deg(U)=10, g(U)=8 \\ \chi({\mathcal O}_U)= -4 \\ \operatorname{gens\,id.}: {4}^{2}\ {5}^{5}\  \end{array} $  &  $\begin{array}{ll} \operatorname{num\,lines}: & 4 \\ \deg\operatorname{others}: & 0 \end{array}$  &  $\begin{array}{c} g(\widetilde{U}) = 8 \\ \deg(\widetilde{U})= 14 \end{array} $  &  $\det \left(\!\begin{array}{cc}
14&9\\
9&2
\end{array}\!\right) = -53 $ & 7\\
\hline
8 & $\mathfrak{S}\left(\!\begin{array}{cccc}
5&7&2&0\\
2&4&1&0
\end{array}\!\right) $ & $\begin{array}{c} \deg(S)=10 , g(S)= 4 \\ \operatorname{gens\,id.}: {3}^{7}\ {4}^{2}\  \\ \deg(C)= 4 , n(C)= 3 \end{array}$  & $53 = 6 \cdot 8 + 5 =  \det \left(\!\begin{array}{ccc}
3&10&1\\
10&42&1\\
1&1&3
\end{array}\!\right)$  &  $\mathcal C_{26} \cap \mathcal C_8 $  &  $3 = 54 - (51 + 6 - 6)$ & $51,  - $ & $7,7$ & $6$ &  -  & $ d_1(\mu) = 4 $ &  $\begin{array}{c} W\subset\mathbb{P}^{9} \\ \deg(W)=9, g(W)=3 \\ \operatorname{gens\,id.}: {2}^{12}\  \end{array} $  &  $\begin{array}{c} \deg(U)=25, g(U)=17 \\ \chi({\mathcal O}_U)= -1 \\ \operatorname{gens\,id.}: {2}^{13}\ {3}^{14}\  \end{array} $  &  $\begin{array}{ll} \operatorname{num\,lines}: & 5 \\ \deg\operatorname{others}: & 2 \end{array}$  &  $\begin{array}{c} g(\widetilde{U}) = 18 \\ \deg(\widetilde{U})= 34 \end{array} $  &  $\det \left(\!\begin{array}{cc}
34&11\\
11&2
\end{array}\!\right) = -53 $ & 8\\
\hline
9 & $\mathfrak{S}\left(\!\begin{array}{cccc}
5&6&2&0\\
2&3&1&0
\end{array}\!\right) $ & $\begin{array}{c} \deg(S)=11 , g(S)= 4 \\ \operatorname{gens\,id.}: {3}^{4}\ {4}^{9}\ {5}^{1}\  \\ \deg(C)= 5 , n(C)= 6 \end{array}$  & $61 = 7 \cdot 8 + 5 =  \det \left(\!\begin{array}{ccc}
3&11&1\\
11&53&0\\
1&0&3
\end{array}\!\right)$  &  $\mathcal C_{38} \cap \mathcal C_8 $  &  $2 = 54 - (53 + 3 - 4)$ & $53, 55$ & $4,4$ & $4$ & $1$ & $ d_1(\mu) = 4 $ &  $\begin{array}{c} W\subset\mathbb{P}^{8} \\ \deg(W)=8, g(W)=3 \\ \operatorname{gens\,id.}: {2}^{7}\  \end{array} $  &  $\begin{array}{c} \deg(U)=22, g(U)=15 \\ \chi({\mathcal O}_U)= -1 \\ \operatorname{gens\,id.}: {2}^{8}\ {3}^{12}\  \end{array} $  &  $\begin{array}{ll} \operatorname{num\,lines}: & 4 \\ \deg\operatorname{others}: & 2 \end{array}$  &  $\begin{array}{c} g(\widetilde{U}) = 16 \\ \deg(\widetilde{U})= 30 \end{array} $  &  $\det \left(\!\begin{array}{cc}
30&11\\
11&2
\end{array}\!\right) = -61 $ & 9\\
\hline
10 & $\mathfrak{S}\left(\!\begin{array}{cccc}
6&2&6&0\\
1&2&0&0
\end{array}\!\right) $ & $\begin{array}{c} \deg(S)=10 , g(S)= 4 \\ \operatorname{gens\,id.}: {3}^{7}\ {4}^{2}\  \\ \deg(C)= 4 , n(C)= 3 \end{array}$  & $69 = 8 \cdot 8 + 5 =  \det \left(\!\begin{array}{ccc}
3&10&1\\
10&44&1\\
1&1&3
\end{array}\!\right)$  &  $\mathcal C_{32} \cap \mathcal C_8 $  &  $3 = 54 - (49 + 6 - 4)$ & $49, 52$ & $7,7$ & $4$ & $1$ & $ d_1(\mu) = 4 $ &  $\begin{array}{c} W\subset\mathbb{P}^{9} \\ \deg(W)=9, g(W)=3 \\ \operatorname{gens\,id.}: {2}^{12}\  \end{array} $  &  $\begin{array}{c} \deg(U)=24, g(U)=14 \\ \chi({\mathcal O}_U)= -4 \\ \operatorname{gens\,id.}: {2}^{13}\ {3}^{14}\  \end{array} $  &  $\begin{array}{ll} \operatorname{num\,lines}: & 2 \\ \deg\operatorname{others}: & 0 \end{array}$  &  $\begin{array}{c} g(\widetilde{U}) = 14 \\ \deg(\widetilde{U})= 26 \end{array} $  &  $\det \left(\!\begin{array}{cc}
26&11\\
11&2
\end{array}\!\right) = -69 $ & 10\\
\hline
11 & $\mathfrak{S}\left(\!\begin{array}{cccc}
6&4&5&0\\
2&3&2&0
\end{array}\!\right) $ & $\begin{array}{c} \deg(S)=12 , g(S)= 5 \\ \operatorname{gens\,id.}: {3}^{1}\ {4}^{21}\ {5}^{1}\  \\ \deg(C)= 5 , n(C)= 6 \end{array}$  & $69 = 8 \cdot 8 + 5 =  \det \left(\!\begin{array}{ccc}
3&12&1\\
12&60&1\\
1&1&3
\end{array}\!\right)$  &  $\mathcal C_{36} \cap \mathcal C_8 $  &  $3 = 54 - (54 + 0 - 3)$ & $54, 57$ & $1,1$ & $3$ & $1$ & $ d_1(\mu) = 4 $ &  $\begin{array}{c} W\subset\mathbb{P}^{5} \\ \deg(W)=3, g(W)=1 \\ \operatorname{gens\,id.}: {3}^{1}\  \end{array} $  &  $\begin{array}{c} \deg(U)=15, g(U)=12 \\ \chi({\mathcal O}_U)= -2 \\ \operatorname{gens\,id.}: {3}^{1}\ {4}^{16}\  \end{array} $  &  $\begin{array}{ll} \operatorname{num\,lines}: & 3 \\ \deg\operatorname{others}: & 4 \end{array}$  &  $\begin{array}{c} g(\widetilde{U}) = 14 \\ \deg(\widetilde{U})= 26 \end{array} $  &  $\det \left(\!\begin{array}{cc}
26&11\\
11&2
\end{array}\!\right) = -69 $ & 11\\
\hline
12 & $\mathfrak{S}\left(\!\begin{array}{cccc}
6&1&6&0\\
1&1&0&0
\end{array}\!\right) $ & $\begin{array}{c} \deg(S)=11 , g(S)= 4 \\ \operatorname{gens\,id.}: {3}^{4}\ {4}^{9}\ {5}^{1}\  \\ \deg(C)= 5 , n(C)= 6 \end{array}$  & $77 = 9 \cdot 8 + 5 =  \det \left(\!\begin{array}{ccc}
3&11&1\\
11&55&0\\
1&0&3
\end{array}\!\right)$  &  $\mathcal C_{44} \cap \mathcal C_8 $  &  $2 = 54 - (51 + 3 - 2)$ & $51, 53$ & $4,4$ & $2$ & $1$ & $ d_1(\mu) = 4 $ &  $\begin{array}{c} W\subset\mathbb{P}^{8} \\ \deg(W)=8, g(W)=3 \\ \operatorname{gens\,id.}: {2}^{7}\  \end{array} $  &  $\begin{array}{c} \deg(U)=21, g(U)=12 \\ \chi({\mathcal O}_U)= -4 \\ \operatorname{gens\,id.}: {2}^{8}\ {3}^{12}\  \end{array} $  &  $\begin{array}{ll} \operatorname{num\,lines}: & 1 \\ \deg\operatorname{others}: & 0 \end{array}$  &  $\begin{array}{c} g(\widetilde{U}) = 12 \\ \deg(\widetilde{U})= 22 \end{array} $  &  $\det \left(\!\begin{array}{cc}
22&11\\
11&2
\end{array}\!\right) = -77 $ & 12\\
\hline
13 & $\mathfrak{S}\left(\!\begin{array}{cccc}
8&0&4&4\\
2&0&4&1
\end{array}\!\right) $ & $\begin{array}{c} \deg(S)=12 , g(S)= 5 \\ \operatorname{gens\,id.}: {3}^{1}\ {4}^{21}\ {5}^{1}\  \\ \deg(C)= 5 , n(C)= 6 \end{array}$  & $85 = 10 \cdot 8 + 5 =  \det \left(\!\begin{array}{ccc}
3&12&1\\
12&62&1\\
1&1&3
\end{array}\!\right)$  &  $\mathcal C_{42} \cap \mathcal C_8 $  &  $3 = 54 - (52 + 0 - 1)$ & $52, 55$ & $1,1$ & $1$ & $1$ & $ d_1(\mu) = 4 $ &  $\begin{array}{c} W\subset\mathbb{P}^{5} \\ \deg(W)=3, g(W)=1 \\ \operatorname{gens\,id.}: {3}^{1}\  \end{array} $  &  $\begin{array}{c} \deg(U)=14, g(U)=9 \\ \chi({\mathcal O}_U)= -5 \\ \operatorname{gens\,id.}: {3}^{1}\ {4}^{17}\  \end{array} $  &  $\begin{array}{ll} \operatorname{num\,lines}: & 0 \\ \deg\operatorname{others}: & 2 \end{array}$  &  - & - & 13
\\ \hline
\end{tabular}
\end{adjustbox}}
\restoregeometry
\clearpage
\pagestyle{headings}

\section*{Acknowledgements}

We thank Sandro Verra for useful discussions. During the revision of the manuscript, we used Claude (Anthropic), an AI assistant, to check computations and proofs and to suggest corrections and improvements to the exposition. We reviewed all the resulting changes, and we take full responsibility for the content of the paper.

\bibliographystyle{amsalpha}
\bibliography{references.bib}

\end{document}